\documentclass[12pt]{article}
\usepackage[english]{babel}
\usepackage{amssymb, amsmath, amsfonts}
\usepackage{mathrsfs}
\newtheorem{Df}{Definition}
\newtheorem{Th}{Theorem}
\newtheorem{Lm}{Lemma}

\newtheorem{Cr}{Corollary}

\begin{document}
\begin{center}
{\bf THE HITTING TIMES WITH TABOO\\
FOR A RANDOM WALK ON AN INTEGER LATTICE}
\end{center}
\vskip0,5cm
\begin{center}
Ekaterina Vl. Bulinskaya\footnote{The work is partially supported by
RFBR grant 10-01-00266.}$^,$\footnote{Lomonosov Moscow State
University.}
\end{center}
\vskip1cm

\begin{abstract}
For a symmetric, homogeneous and irreducible random walk on $\mathbb{Z}^{d}$, $d\in\mathbb{N}$,
having zero mean and a finite variance of jumps, we study the passage times (taking values in $[0,\infty]$)
determined by the starting point $x$, the hitting state $y$ and the taboo state $z$.
We find the probability that these passages times are finite and analyze the tails of
their cumulative distribution functions.
In particular, it turns out that for the random walk on $\mathbb{Z}^{d}$,
except for a simple (nearest neighbor) random walk on $\mathbb{Z}$,
the order of the tail decrease is specified by dimension $d$ only.
In contrast, for a simple random walk on $\mathbb{Z}$, the asymptotic properties of
hitting times with taboo essentially depend on the mutual location of the points $x$, $y$ and $z$.
These problems originated in our recent study of branching
random walk on $\mathbb{Z}^{d}$ with a single source of branching.

\vskip0,5cm {\it Keywords and phrases}: random walks on integer
lattices, hitting times, taboo probabilities, branching random walk. \vskip0,5cm
{\it $2010$ AMS classification}: 60G50, 60J27, 60G17.
\end{abstract}

\section{Introduction}
A random walk is a classical model being the source of
numerous interesting problems having elegant solutions.
The monographs \cite{Lawler}, \cite{Madras_Slade} and \cite{Spitzer}
devoted to various properties of random walks have become the reference books for many researchers.
The recent publications such as \cite{Chen} and \cite{Lawler_Limic} demonstrate
its non-vanishing popularity.
The reason is the abundance of applications of this simple probabilistic model (see, e.g.,
\cite{Berg} and \cite{Weiss}). Moreover,
a lot of new complicated models were constructed on the basis of random walk
(one can mention those investigated in \cite{Caravenna} and \cite{Revesz}).
Thus, many assertions
of the theory of random walks have been obtained as auxiliary (technical) lemmas needed for understanding other
models of interest.

The present paper also provides results concerning random walk on integer lattice $\mathbb{Z}^{d}$
($d\in\mathbb{N}$)
which arose during the study of another model, namely, branching random walk (BRW) on
$\mathbb{Z}^{d}$ with
a single source of branching (different modifications of the model were considered, e.g.,
in papers \cite{Albeverio_et_al}, \cite{MZ_Bulinskaya}, \cite{VT_Siberia} and \cite{Y_TVP}).
Its main features are the following.
If the starting point of the process is $x\in\mathbb{Z}^{d}$, $x\neq{\bf 0}$, then
the parent particle performs random walk on $\mathbb{Z}^{d}$ until hitting the origin.
If the starting point is ${\bf 0}$, or the parent particle has just hit it, then after exponentially
distributed time the particle may die producing a random number of offsprings.
Otherwise, it leaves the origin (jumps randomly to a point $x'\in\mathbb{Z}^{d}\backslash\{{\bf 0}\}$)
and behaves afterwards as a memoryless particle starting at $x'$.
At the birth moment the newborn particles are located at the origin.
They and all their descendants evolve according to the scheme described above
independently of each other as well as of their parents history.

In the framework of BRW a natural question arises:
what is the limit distribution, as $t\to\infty$, of the number of particles $\mu(t;y)$
at each point $y\in\mathbb{Z}^{d}$ at time $t$?
Evidently, $\mu(t;y)$ essentially depends on
the trajectory of the parent particle. Indeed, if the parent particle hits the origin
earlier than the time $t$ then $\mu(t;y)$ can take any value in $\mathbb{Z}_{+}$.
Otherwise, $\mu(t;y)=1$ or $\mu(t;y)=0$ if the parent particle is located at time $t$ at point $y$
or outside it, respectively.

This is the motivation to introduce
the notion of hitting time with taboo for a random walk on $\mathbb{Z}^{d}$.
More precisely, for any $x,y,z\in\mathbb{Z}^{d}$ such that $y\neq z$, let function $H_{x,y,z}(t),$
$t\geq0,$ be (improper) cumulative distribution function (c.d.f.) of the first hitting time (or the first
return time if $x=y$)
of point $y$ if the starting point of the random walk is $x$ and the point $z$ is a taboo state.
In the next section we give a formal definition of this function. In the present work we find
the limit value $H_{x,y,z}(\infty)=\lim\nolimits_{t\to\infty}{H_{x,y,z}(t)}$ and analyze
the asymptotic behavior of $H_{x,y,z}(\infty)-H_{x,y,z}(t)$, as $t\to\infty$,
for a symmetric, homogeneous, irreducible random walk on $\mathbb{Z}^{d}$
($d\in\mathbb{N}$) having zero mean and a finite variance of jumps.

It turns out that the most interesting case is $d=1$ since there are
two quite different kinds of asymptotic behavior of $H_{x,y,z}(t)$, $t\to\infty$,
depending on whether the random walk is simple (nearest neighbor) or not.
It is also worth mentioning that, for the random walk on $\mathbb{Z}^{d}$,
except for a simple random walk on $\mathbb{Z}$, the value $H_{x,y,z}(\infty)\in(0,1)$
and the order of decrease of function $H_{x,y,z}(\infty)-H_{x,y,z}(t)$ is
determined by dimension $d$ only, regardless of $x,y$ and $z$. In contrast,
for a simple random walk on $\mathbb{Z}$,
the mutual location of points $x,y$ and $z$ determines
the value $H_{x,y,z}(\infty)\in[0,1]$
as well as the order of decrease of $H_{x,y,z}(\infty)-H_{x,y,z}(t)$ as $t\to\infty$.

Finally, we recall that the properties of hitting times
(or, more generally, passage times including first entrance and last exit times)
and taboo probabilities for a Markov chain with stationary transition
probabilities were exposed in \cite{Chung} and \cite{Syski}, respectively.
For the taboo probabilities see also \cite{Zubkov} and references therein.
The counterpart of the function $H_{x,y,z}(t)$ for a Markov chain with stationary transition
probabilities has been used earlier
(see, e.g., \cite{Bulinski_Shiryaev}, p.~202, and \cite{Karlin_Taylor}, p.~31).
However, it was treated as an auxiliary tool and has not been studied per se.

\section{Main results}

We assume that all random variables (taking values in $[-\infty,\infty]$) are defined on
a probability space $(\Omega,\mathcal{F},{\sf P})$.
We study a class of random walks on $\mathbb{Z}^{d}$ more general than that
proposed in \cite{Lawler_Limic}, Ch.1, Sec.2.

\begin{Df}
For $d\in\mathbb{N}$, let $A=(a(x,y))_{x,y\in\mathbb{Z}^{d}}$ be a matrix such that elements
$a(x,y)\geq0$ whenever $x\neq y$, $a(x,x)<0$, $\sum\nolimits_{y\in\mathbb{Z}^{d}}{a(x,y)}=0$
and $\sup\nolimits_{x\in\mathbb{Z}^{d}}{|a(x,x)|}<\infty$.
Then a continuous time Markov chain $S=\{S(t),t\geq0\}$
having the state space $\mathbb{Z}^{d}$, i.e. $S(t):\Omega\rightarrow\mathbb{Z}^{d}$ for each $t\geq0$, and
generator $A$ is called a random walk on $\mathbb{Z}^{d}$.
\end{Df}
The existence of such a Markov chain is guarantied by Theorem 2 in \cite{Gikhman_Skorokhod}, Ch.3, Sec.2.

From here on we consider only a symmetric, homogeneous, irreducible random walk $S$
having a finite variance of jumps. \emph{Symmetry} and \emph{homogeneity}  mean that
for $x,y\in\mathbb{Z}^{d}$ one has ${a(x,y)=a(y,x)}$ and ${a(x,y)=a({\bf 0},y-x)=:a(y-x)}$, respectively.
Random walk is called \emph{irreducible} if for any $x,y\in\mathbb{Z}^{d}$ there exists
$t>0$ such that ${\sf P}(S(t)=y|S(0)=x)>0$, i.e. all the points of the lattice
$\mathbb{Z}^{d}$ can be reached. Furthermore, a random walk is said to have
a \emph{finite variance of jumps} if $\sum\nolimits_{x\in\mathbb{Z}^{d}}{\|x\|^{2}a({\bf 0},x)<\infty}$
where ${\bf 0}$ is the origin of $\mathbb{Z}^{d}$ and $\|\cdot\|$ is a norm in this space.

The conditions imposed on elements $a(x,y)$ allow us to use an explicit construction of the random walk
on $\mathbb{Z}^{d}$ with generator $A$ (see, e.g., Theorem 1.2 in \cite{Bremaud}, Ch.9, Sec.1).
According to this construction $S$ is a regular jump process with right continuous trajectories.
A jump of the process from state $x\in\mathbb{Z}^{d}$ to state $y\in\mathbb{Z}^{d}$, $x\neq y$,
occurs with probability $a(x,y)/a(x,x)$ independently of the earlier evolution of $S$. Moreover,
for transition times of the process $\tau^{(0)}:=0$ and
$\tau^{(n)}:=\inf\{t\geq\tau^{(n-1)}:S(t)\neq S(\tau^{(n-1)})\}$, $n\geq1$, the following statement holds.
Random variables $\{\tau^{(n+1)}-\tau^{(n)}\}_{n=0}^{\infty}$
are independent and each of them has exponential distribution with parameter $a:=-a({\bf 0},{\bf 0})$.
In what follows we consider the version of the process $S$ constructed in this way.

Note that for the sequence $\{\tau^{(n)}\}_{n=1}^{\infty}$ of transition times of $S$
one can consider {\it embedded chain} $\{S_{n},n\in\mathbb{Z}_{+}\}$ defined by way of $S_{n}:=S(\tau^{(n)})$.
Mention in passing that ${\{S_{n},n\in\mathbb{Z}_{+}\}}$ is a discrete-time homogeneous Markov chain
having transition probabilities $a(x,y)/a$ for $x\neq y$, $x,y\in\mathbb{Z}^{d}$
(see, e.g., \cite{Bremaud}, Ch.8, Sec.4).

Recall some other facts concerning the introduced random walk on $\mathbb{Z}^{d}$.
Let $p(t;x,y),$ $t\geq0$, $x,y\in\mathbb{Z}^{d},$ be the \emph{transition probabilities} of
the random walk, that is,
$$p(t;x,y):={\sf P}(S(t+u)=y|S(u)=x)={\sf P}(S(t)=y|S(0)=x)$$
for any $u\geq0$. Then for $h\to0+$
\begin{eqnarray*}
p(h;x,y)&=&a(x,y)h+o(h),\quad\mbox{if}\quad x\neq y,\\
p(h;x,x)&=&1+a(x,x)h+o(h).
\end{eqnarray*}
Due to \cite{Bremaud}, Ch.8, Sec.3,
the transition probabilities of the process $S$ satisfy the backward Kolmogorov equations
\begin{equation}\label{backward_K_equations}
\frac{d p(t;x,y)}{d t}=\sum\nolimits_{z\in\mathbb{Z}^{d}}{a(x,z) p(t;z,y)},\quad p(0;x,y)=\delta_{y}(x).
\end{equation}
As usual, $\delta_{y}(x)$ is equal to $1$ or $0$ for $x=y$ or $x\neq y$, respectively.
The Fourier transform applied to \eqref{backward_K_equations}
yields (see \cite{VT_Siberia} and \cite{Yarovaya_book}, Ch.2, Sec.1)
\begin{equation}\label{p(t;x,y)=int}
p(t;x,y)=\frac{1}{(2\pi)^{d}}\int\nolimits_{[-\pi,\pi]^{d}}{e^{\phi(\theta) t+i(\theta,y-x)}\,d\theta},
\quad x,y\in\mathbb{Z}^{d},
\end{equation}
where $(\cdot,\cdot)$ stands for the scalar product in Euclidean space $\mathbb{R}^{d}$ and
$$\phi(\theta):=\sum\limits_{z\in\mathbb{Z}^{d}}{a({\bf 0},z)\cos(z,\theta)},\quad\theta\in[-\pi,\pi]^{d}.$$
Furthermore, by the method of steepest descent one derives from \eqref{p(t;x,y)=int} the asymptotic relation
\begin{equation}\label{p(t;x,y)sim}
p(t;x,y)\sim\frac{\gamma_{d}}{t^{d/2}},\quad x,y\in\mathbb{Z}^{d},\quad t\to\infty.
\end{equation}
Here
$$\gamma_{d}:=\frac{1}{(2\pi)^{d/2}\sqrt{|\det{\phi''_{\theta \theta}({\bf 0})}|}}\quad\mbox{and}\quad
\phi''_{\theta \theta}({\bf 0}):=\left(\left.\frac{\partial^{2}\phi(\theta)}{\partial \theta_{i} \partial
\theta_{j}}\right|_{\theta={\bf 0}}\right)_{i,j\in\{1,\ldots,d\}}.$$
A similar result for transition probabilities of \emph{a discrete time} random walk on $\mathbb{Z}^{d}$
was established in \cite{Spitzer}, Ch.2, Sec.7.

Moreover, according to \cite{Yarovaya_book}, Ch.2, Sec.1, equality \eqref{p(t;x,y)=int} implies
that
\begin{equation}\label{p(t;0,0)-p(t,x,y)sim}
p(t;{\bf 0},{\bf 0})-p(t;x,y)\sim\frac{\tilde{\gamma}_{d}(y-x)}{t^{d/2+1}},
\quad x,y\in\mathbb{Z}^{d},\quad t\to\infty,
\end{equation}
where
$$\tilde{\gamma}_{d}(z):=\frac{1}{2(2\pi)^{d}}
\int\nolimits_{\mathbb{R}^{d}}{(\upsilon,z)^{2}e^{\frac{1}{2}(\phi''_{\theta \theta}({\bf 0})\upsilon,\upsilon)}
\,d\upsilon},\quad z\in\mathbb{Z}^{d}.$$

Set $G_{\lambda}(x,y):=\int\nolimits_{0}^{\infty}{e^{-\lambda t}p(t;x,y)\,d t},$ $\lambda\geq0$,
$x,y\in\mathbb{Z}^{d}$, that is, $G_{\lambda}(x,y)$ is the Laplace transform of the transition
probability $p(\cdot;x,y)$. In view of \eqref{p(t;x,y)=int} the introduced function
can be represented in the form
\begin{equation}\label{G_lambda(x,y)=}
G_{\lambda}(x,y)=\frac{1}{(2\pi)^{d}}\int\nolimits_{[-\pi,\pi]^{d}}{\frac{e^{i(\theta,y-x)}}
{\lambda-\phi(\theta)}\,d\theta}.
\end{equation}

Taking into account \eqref{p(t;x,y)sim} it is easily seen that for $d\geq3$ the \emph{Green's function}
$G_{0}(x,y)$ is finite for all $x,y\in\mathbb{Z}^{d}$. However, for $d=1$ or $d=2$, one has
$\lim\nolimits_{\lambda\to0+}{G_{\lambda}(x,y)}=\infty$. This phenomenon corresponds
to transience of the random walk for $d\geq3$ and its recurrence for $d=1$ or $d=2$
(for more details see, e.g., \cite{Lawler_Limic}, Ch.4, Sec.1 and 2).

By virtue of \eqref{p(t;0,0)-p(t,x,y)sim}
the function $\lim\nolimits_{\lambda\to0+}{(G_{\lambda}({\bf 0},{\bf 0})-G_{\lambda}(x,y))}$ is finite
for all $d\in\mathbb{N}$ and $x,y\in\mathbb{Z}^{d}$.
Therefore we can define
$$\rho_{d}(x):=\left\{
\begin{array}{lcl}
a\lim\limits_{\lambda\to0+}(G_{\lambda}({\bf 0},{\bf
0})-G_{\lambda}({\bf 0},x)),&\mbox{if}&\quad
x\neq{\bf 0},\\
1,&\mbox{if}&\quad x={\bf 0}.
\end{array}
\right.
$$
Equality \eqref{G_lambda(x,y)=} allows us to give another formula for $\rho_{d}(x)$ if $x\neq{\bf 0}$
\begin{equation}\label{rho_d(x)=int}
\rho_{d}(x)=\frac{a}{(2\pi)^{d}}\int\nolimits_{[-\pi,\pi]^{d}}{\,\frac{\cos(x,\theta)-1}{\phi(\theta)}\,d\theta}.
\end{equation}

Now we are able to introduce some basic notation and formulate our main results.
Set $\tau:=\inf\{t\geq0: S(t)\neq S(0)\}$, i.e. the stopping time $\tau$
(with respect to the natural filtration of the process $S$) is the time of the first
exit from the starting point of the random walk. Taking into account the explicit construction of the random walk
one has $\tau=\tau^{(1)}$ almost surely (a.s.) and, consequently, $G(t):={\sf P}(\tau\leq t|S(0)=x)=1-e^{-a t}$, $t\geq0$.

Let $\tau_{y,z}:=\inf\{t\geq\tau: S(t)=y, S(u)\neq z, \tau\leq u\leq t\}$ for $y,z\in\mathbb{Z}^{d}$, $y\neq z$. As usual,
$\inf\{t\in\varnothing\}=\infty$.
We call the stopping time $\tau_{y,z}$ a \emph{hitting time of the state} $y$ \emph{with the taboo state} $z$.
Denote its (improper) c.d.f. given the starting point $x$
by $H_{x,y,z}(t):={\sf P}(\tau_{y,z}\leq t|S(0)=x)$, $x,y,z\in\mathbb{Z}^{d}$, $y\neq z$, $t\geq0$.

Note that the random variable $\tau^{-}_{y,z}:=\inf\{t\geq0:S(t+\tau)=y, S(u)\neq z,
\tau\leq u\leq t+\tau\}$, $y,z\in\mathbb{Z}^{d}$, $y\neq z$,
is not a stopping time w.r.t. the natural filtration of the process $S$.
However, $\tau_{y,z}$ can be also called a
\emph{hitting time of the state} $y$ \emph{with the taboo state} $z$, since
in this case we count time not from the moment of the process start but from the moment
of the first exit out of the starting point. In a similar way, we set
$H^{-}_{x,y,z}(t):={\sf P}(\tau^{-}_{y,z}\leq t|S(0)=x),$ $x,y,z\in\mathbb{Z}^{d}$, $y\neq z$, $t\geq0$.

Obviously, $\tau_{y,z}=\tau^{-}_{y,z}+\tau$ a.s.
Moreover, in view of strong Markov property of the random walk w.r.t. the stopping time $\tau$
(see, e.g., \cite{Bremaud}, Ch.8, Sec.4) the random variables $\tau^{-}_{y,z}$ and $\tau$
are independent. Hence, $\tau_{y,z}$ has an absolutely continuous (improper)
distribution as a sum of two independent random variables, one of them having a density.
On the other hand, the explicit construction of $S$ shows that the distribution of $\tau^{-}_{y,z}$
has an atom at zero when $x\neq y$, namely, $H^{-}_{x,y,z}(0)=a(x,y)/a$.

For $x,y\in\mathbb{Z}^{d}$ such that $y\neq{\bf 0}$, we specify the functions
$C_{d}(x,y)$ by formulae
$$C_{1}(x,y):=\frac{\rho_{1}(y-x)+\rho_{1}(x)-\rho_{1}(y)}{4a\pi\gamma_{1}}
+\frac{a \pi y^{2}\gamma^{3}_{1}(\rho_{1}(y-x)-\rho_{1}(x)-\rho_{1}(y))}{\rho^{2}_{1}(y)}
+\frac{2a \pi xy\gamma^{3}_{1}}{\rho_{1}(y)},$$
$$C_{2}(x,y):=\frac{\rho_{2}(y-x)+\rho_{2}(x)-\rho_{2}(y)}{4 a \gamma_{2}},$$
$$C_{d}(x,y):=\frac{2\gamma_{d}(\rho_{d}(y-x)+\rho_{d}(x)-\rho_{d}(y))}{a(d-2)(G_{0}({\bf 0},{\bf 0})
+G_{0}({\bf 0},y))^{2}},\quad d\geq3.$$

The following theorems are the main results of the paper.

\begin{Th}\label{T:not_simple}
Let $x,y,z\in\mathbb{Z}^{d}$ be such that $y\neq z$.
Then for the random walk $S$ on $\mathbb{Z}^{d}$, except for a simple random walk on $\mathbb{Z}$, one has
\begin{eqnarray}
H_{x,y,z}(\infty)&=&\frac{\rho_{d}(x-z)+\rho_{d}(y-z)-\rho_{d}(y-x)}{2\rho_{d}(y-z)}\in(0,1),
\quad d=1\quad\mbox{or}\quad d=2,\label{T_Hxyz(infty),d=1,2}\\
H_{x,y,z}(\infty)&=&\frac{G_{0}({\bf 0},{\bf 0})\rho_{d}(y-z)-G_{0}({\bf 0},{\bf 0})\rho_{d}(y-x)
+G_{0}(y,z)\rho_{d}(x-z)}{\rho_{d}(y-z)(G_{0}({\bf 0},{\bf 0})+G_{0}(y,z))}\in(0,1),\; d\geq3,\quad\;
\label{T_Hxyz(infty),d>=3}
\end{eqnarray}
Moreover, as $t\to\infty$,
\begin{eqnarray}
H_{x,y,z}(\infty)-H_{x,y,z}(t)&\sim&\frac{C_{1}(x-z,y-z)}{\sqrt{t}}\quad\mbox{for}\quad d=1,
\label{H_xyz(infty)-H_xyz(t)sim,d=1}\\
H_{x,y,z}(\infty)-H_{x,y,z}(t)&\sim&\frac{C_{2}(x-z,y-z)}{\ln{t}}\quad\mbox{for}\quad d=2,
\label{H_xyz(infty)-H_xyz(t)sim,d=2}\\
H_{x,y,z}(\infty)-H_{x,y,z}(t)&\sim&\frac{C_{d}(x-z,y-z)}{t^{d/2-1}}\quad\mbox{for}\quad d\geq3
\label{H_xyz(infty)-H_xyz(t)sim,d>=3}
\end{eqnarray}
where $C_{d}(\cdot,\cdot)$, $d\in\mathbb{N}$, are positive functions defined above.
\end{Th}

\begin{Th}\label{T:simple}
Let $S$ be a simple random walk on $\mathbb{Z}$. If $x,y,z\in\mathbb{Z}$ and $y\neq z$ then
\begin{eqnarray}
1-H_{x,y,z}(t)\sim\frac{\sqrt{2}\,|y-x|}{\sqrt{a\,\pi}\,\sqrt{t}},& & x<y<z\quad
\mbox{or}\quad z<y<x,\label{Hxyz_simple_1}\\
1-\frac{1}{2|y-z|}-H_{y,y,z}(t)\sim\frac{1}{\sqrt{2\,a\,\pi}\,\sqrt{t}},& & x=y,
\label{Hxyz_simple_2}\\
\frac{x-z}{y-z}-H_{x,y,z}(t)=o(e^{-a\,\varepsilon\,t}),& &z<x<y\quad\mbox{or}
\quad y<x<z,\label{Hxyz_simple_3}\\
\frac{1}{2|y-z|}-H_{z,y,z}(t)=o(e^{-a\,\varepsilon\,t}),& &x=z,
\label{Hxyz_simple_4}\\
H_{x,y,z}(t)\equiv0,& &x<z<y\quad\mbox{or}\quad y<z<x,\label{Hxyz_simple_5}
\end{eqnarray}
for some $\varepsilon\in(0,1)$.
\end{Th}

\begin{Th}\label{T:H-}
Theorems \emph{\ref{T:not_simple}} and \emph{\ref{T:simple}} hold true if function $H_{x,y,z}$
at the left-hand sides of $\eqref{T_Hxyz(infty),d=1,2}$--$\eqref{Hxyz_simple_5}$
is replaced by $H^{-}_{x,y,z}$, whereas the right-hand sides of these formulae remain intact.
\end{Th}

The proofs of Theorems \ref{T:not_simple} and \ref{T:simple} are given in Sections 4 and 5, respectively.
They bear on application of the Laplace-Stieltjes
and the Laplace transforms to the main equations involving $H_{x,y,z}$ and $H_{x,z,y}$ as well as functions
$H_{x,y}$, $H_{z,y}$, $H_{x,z}$ and $H_{y,z}$. Here $H_{x,y}(t)$, $x,y\in\mathbb{Z}^{d}$, $t\geq0$,
is the (improper) c.d.f. of the hitting time of point $y$ given the starting point $x$.
After that we employ the Tauberian theorems techniques and, for $d$ sufficiently large, we essentially use
Fa\`{a} di Bruno's formula for the $n$-th derivative of two functions superposition, $n\geq1$.
An essential difficulty arising on this way is to prove the positivity of functions $C_{d}(x,y)$,
$x,y\in\mathbb{Z}^{d}$, $y\neq{\bf 0}$,
under the conditions of Theorem~\ref{T:not_simple}. Theorem \ref{T:H-} proved in Section 6
can be viewed as a corollary of Theorems~\ref{T:not_simple} and \ref{T:simple}.

\section{Auxiliary results}

For a nonnegative function $\kappa(t)$ and a nonnegative nondecreasing function $\chi(t)$, $t\geq0,$
we set
$$\hat{\kappa}(\lambda):=\int\nolimits_{0}^{\infty}{e^{-\lambda t}\kappa(t)\,dt},\quad
\check{\chi}(\lambda):=\int\nolimits_{0}^{\infty}{e^{-\lambda t}\,d \chi(t)},\quad \lambda>0,$$
whenever the integrals exist.
Recall a useful relation linking the Laplace-Stieltjes transform of the bounded function $\chi(t)$
and the Laplace transform of $\chi(\infty)-\chi(t)$ where ${\chi(\infty):=\lim\nolimits_{t\to\infty}{\chi(t)}}$,
namely,
\begin{equation}\label{universal_formula}
\widehat{(\chi(\infty)-\chi)}(\lambda)=\frac{\chi(\infty)-\check{\chi}(\lambda)}{\lambda},\quad\lambda>0.
\end{equation}
This relation results from the formula of integration by parts.

Introduce a stopping time $\tau_{y}:=\inf\{t\geq\tau: S(t)=y\}$ called \emph{the hitting time of point}
$y\in\mathbb{Z}^{d}$. Denote its (conditional) c.d.f.
given the starting point $x$ by
$H_{x,y}(t):={\sf P}(\tau_{y}\leq t|S(0)=x)$, $x,y\in\mathbb{Z}^{d}$, $t\geq0$.
The limit value $H_{x,y}(\infty):=\lim\nolimits_{t\to\infty}{H_{x,y}(t)}$ and
the asymptotic behavior of function $H_{x,y}(\infty)-H_{x,y}(t)$, as $t\to\infty$, were found
in papers \cite{TVP_Bulinskaya}--\cite{LMJ_Bulinskaya}, \cite{VT_Siberia} and \cite{VTY} for $y={\bf 0}$
and $x\in\mathbb{Z}^{d}$. In Lemma \ref{L:Hxy} we reformulate these results in a more general way.

\begin{Lm}\label{L:Hxy}
Let $x,y\in\mathbb{Z}^{d}$. Then
\begin{eqnarray}
H_{x,y}(\infty)=1\quad&\mbox{for}&\quad d=1\quad\mbox{and}\quad d=2,\label{Hxy(infty),d=1,d=2}\\
H_{x,x}(\infty)=1-\frac{1}{a\,G_{0}({\bf 0},{\bf 0})}\quad&\mbox{for}&\quad d\geq3\quad\mbox{and}\quad
y=x,\label{Hxx(infty),d>=3}\\
H_{x,y}(\infty)=\frac{G_{0}(x,y)}{G_{0}({\bf 0},{\bf 0})}\quad&\mbox{for}&\quad d\geq3\quad\mbox{and}\quad
y\neq x.\label{Hxy(infty),d>=3}
\end{eqnarray}
Moreover, as $t\to\infty$, one has
\begin{eqnarray}
1-H_{x,y}(t)\sim\frac{\rho_{1}(y-x)}{a\,\gamma_{1}\,\pi\,\sqrt{t}}\quad&\mbox{for}&\quad d=1,
\label{Hxy(infty)-Hxy(t)sim,d=1}\\
1-H_{x,y}(t)\sim\frac{\rho_{2}(y-x)}{a\,\gamma_{2}\,\ln{t}}\quad&\mbox{for}&\quad d=2,
\label{Hxy(infty)-Hxy(t)sim,d=2}\\
H_{x,y}(\infty)-H_{x,y}(t)\sim\frac{2\gamma_{d}\,\rho_{d}(y-x)}{a\,(d-2)\,
G^{2}_{0}({\bf 0},{\bf 0})\,t^{d/2-1}}\quad&\mbox{for}&\quad d\geq3.
\label{Hxy(infty)-Hxy(t)sim,d>=3}
\end{eqnarray}
\end{Lm}
\textsc{Proof.} In view of definitions of $\tau$ and $\tau_{x}$ as well as
of conditional probability and a Markov chain we see that
\begin{eqnarray*}
p(t;x,x)&=&{\sf P}(S(t)=x,\tau>t|S(0)=x)+{\sf P}(S(t)=x,\tau\leq t|S(0)=x)\\
&=&{\sf P}(\tau>t|S(0)=x)+{\sf P}(S(t)=x,\tau_{x}\leq t|S(0)=x)\\
&=&1-G(t)+\int\nolimits_{0}^{t}{{\sf P}(S(t)=x|\tau_{x}=u,S(0)=x)\,d {\sf P}(\tau_{x}\leq u|S(0)=x)}\\
&=&1-G(t)+\int\nolimits_{0}^{t}{{\sf P}(S(t)=x|S(u)=x)\, d H_{x,x}(u)}.
\end{eqnarray*}
Thus we deduce the following formula
\begin{equation}\label{int_eq_p(t,x,x)}
p(t;x,x)=1-G(t)+\int\nolimits_{0}^{t}{p(t-u;x,x)\,d \,H_{x,x}(u)},\quad x\in\mathbb{Z}^{d},\quad t\geq0.
\end{equation}
Note that \eqref{G_lambda(x,y)=} implies the identity $G_{\lambda}(x,x)=G_{\lambda}({\bf 0},{\bf
0})$ for any $x\in\mathbb{Z}^{d}$. Therefore, application of the Laplace-Stieltjes
transform to \eqref{int_eq_p(t,x,x)} leads to
\begin{equation}\label{lt_Hxx}
\check{H}_{x,x}(\lambda)=1-\frac{\widehat{(1-G)}(\lambda)}{G_{\lambda}(x,x)}=1-\frac{\widehat{(1-G)}(\lambda)}{G_{\lambda}({\bf
0},{\bf 0})}=\check{H}_{{\bf 0},{\bf 0}}(\lambda).
\end{equation}
In particular, by the Laplace-Stieltjes transform uniqueness (see, e.g., \cite{Feller}, Ch.13, Sec.1) it follows that
\begin{equation}\label{Hx,x=H0,0}
H_{x,x}(t)= H_{{\bf 0},{\bf 0}}(t),\quad t\geq0,\quad
x\in\mathbb{Z}^{d}.
\end{equation}
For different $d$, asymptotic properties of $H_{{\bf 0},{\bf 0}}(t)$, $t\to\infty$, were established in
\cite{TVP_Bulinskaya}, \cite{VT_Siberia} and \cite{VTY}. More exactly, the limit value
$H^{-}_{{\bf 0},{\bf 0}}(\infty):=\lim\nolimits_{t\to\infty}{H^{-}_{{\bf 0},{\bf 0}}(t)}$ and
the asymptotic behavior
of $H^{-}_{{\bf 0},{\bf 0}}(\infty)-H^{-}_{{\bf 0},{\bf 0}}(t)$, $t\to\infty$, were found
for the function
$H^{-}_{{\bf 0},{\bf 0}}(t)$, $t\geq0$, related to $H_{{\bf 0},{\bf 0}}$ by
${H_{{\bf 0},{\bf 0}}(t)=G\ast H^{-}_{{\bf 0},{\bf 0}}}(t)$
(as usual, $\ast$ denotes the convolution). However, these results are easily extended to the case of function
$H_{{\bf 0},{\bf 0}}(t)$. Hence, due to \eqref{Hx,x=H0,0} relations \eqref{Hxy(infty),d=1,d=2},
\eqref{Hxx(infty),d>=3} and \eqref{Hxy(infty)-Hxy(t)sim,d=1}--\eqref{Hxy(infty)-Hxy(t)sim,d>=3}
are proved for $x=y$.

Let $x\neq y$. Combined definitions of $\tau_{y}$, conditional probability and
a Markov chain imply
\begin{eqnarray*}
p(t;x,y)&=&{\sf P}(S(t)=y,\tau_{y}\leq t|S(0)=x)\\
&=&\int\nolimits_{0}^{t}{{\sf P}(S(t)=y|\tau_{y}=u,S(0)=x)\,d{\sf P}(\tau_{y}\leq u|S(0)=x)}\\
&=&\int\nolimits_{0}^{t}{{\sf P}(S(t)=y|S(u)=y)\,d H_{x,y}(u)}.
\end{eqnarray*}
Consequently, for $x\neq y$ and $t\geq0$ we get the formula
\begin{equation}\label{int_eq_p(t,x,y)}
p(t;x,y)=\int\nolimits_{0}^{t}{p(t-u;y,y)\,d\,H_{x,y}(u)}.
\end{equation}
By virtue of \eqref{G_lambda(x,y)=} identities $G_{\lambda}(x,y)=G_{\lambda}(x-y,{\bf 0})=G_{\lambda}(y-x,{\bf 0})$
and $G_{\lambda}(y,y)=G_{\lambda}({\bf 0},{\bf 0})$ are valid. Then application of the
Laplace-Stieltjes transform
to (\ref{int_eq_p(t,x,y)}) yields
\begin{equation}\label{lt_Hxy}
\check{H}_{x,y}(\lambda)=\frac{G_{\lambda}(x,y)}{G_{\lambda}(y,y)}=\frac{G_{\lambda}(x-y,{\bf
0})}{G_{\lambda}({\bf 0},{\bf 0})}=\frac{G_{\lambda}(y-x,{\bf
0})}{G_{\lambda}({\bf 0},{\bf 0})}=\check{H}_{x-y,{\bf
0}}(\lambda)=\check{H}_{y-x,{\bf 0}}(\lambda).
\end{equation}
Now using the Laplace-Stieltjes transform uniqueness we obtain
\begin{equation}\label{Hx,y=Hx-y,0=Hy-x,0}
H_{x,y}(t)=H_{x-y,{\bf 0}}(t)=H_{y-x,{\bf 0}}(t),\quad
t\geq0.
\end{equation}
The limit value $H_{z,{\bf 0}}(\infty)$, $z\in\mathbb{Z}^{d}\backslash\{{\bf 0}\}$, and the asymptotic
behavior of $H_{z,{\bf 0}}(\infty)-H_{z,{\bf 0}}(t)$, $t\to\infty$, were found in
\cite{TSP_B} and \cite{LMJ_Bulinskaya}. It follows that by \eqref{Hx,y=Hx-y,0=Hy-x,0}
relations \eqref{Hxy(infty),d=1,d=2}, \eqref{Hxy(infty),d>=3} and
\eqref{Hxy(infty)-Hxy(t)sim,d=1}--\eqref{Hxy(infty)-Hxy(t)sim,d>=3} are established for $x\neq y$.
The proof of Lemma \ref{L:Hxy} is completed.

The next lemma provides a linear integral equation involving functions $H_{x,y,z}(t)$ and $H_{x,z,y}(t)$
which is the basis of our theorems proofs.

\begin{Lm}\label{L:main_equation}
Let $x,y,z\in\mathbb{Z}^{d}$ and $y\neq z$. Then the
following equation holds
\begin{equation}\label{main_equation}
H_{x,y}(t)=H_{x,y,z}(t)+\int\nolimits_{0}^{t}{H_{z,y}(t-u)\,d H_{x,z,y}(u)},\quad t\geq0.
\end{equation}
\end{Lm}
\textsc{Proof. }Firstly note that the random variables $\tau_{y}$ and $\tau_{y,z}$ coincide a.s.
on event $\{\tau_{y,z}<\infty\}$ in view of their definition, that is
$\tau_{y}\mathbb{I}(\tau_{y,z}<\infty)=\tau_{y,z}\mathbb{I}(\tau_{y,z}<\infty)$ a.s.
(here $\mathbb{I}(B)$ stands for the indicator of a set $B$).
Consequently, $\{\tau_{y}\leq t,\tau_{y,z}\leq t\}=\{\tau_{y,z}\leq t\}$ for each $t\geq 0$.
Moreover, one can see that $\tau_{y}\mathbb{I}(\tau_{z,y}<\infty)\geq\tau_{z,y}\mathbb{I}(\tau_{z,y}<\infty)$
a.s. Therefore, for every nonnegative $t$, one has
$\{\tau_{y}\leq t, \tau_{y,z}>t\}=\{\tau_{y}\leq t, \tau_{y,z}=\infty\}
=\{\tau_{y}\leq t, \tau_{z,y}<\infty\}=\{\tau_{y}\leq t, \tau_{z,y}\leq t\}$.
These relations allow us to write
\begin{eqnarray*}
H_{x,y}(t)&=&{\sf P}(\tau_{y}\leq t, \tau_{y,z}\leq t|S(0)=x)+{\sf P}(\tau_{y}\leq t, \tau_{y,z}>t|S(0)=x)\\
&=&{\sf P}(\tau_{y,z}\leq t|S(0)=x)+{\sf P}(\tau_{y}\leq t, \tau_{z,y}\leq t|S(0)=x)\\
&=&H_{x,y,z}(t)+\int\nolimits_{0}^{t}{{\sf P}(\tau_{y}\leq t|\tau_{z,y}=u, S(0)=x)\,d H_{x,z,y}(u)}.
\end{eqnarray*}
Using definitions of $\tau_{y}$ and $\tau_{z,y}$ as well as that of a Markov chain
it can be verified that
${{\sf P}(\tau_{y}\leq t|\tau_{z,y}=u,S(0)=x)={\sf P}(\tau_{y}\leq t-u|S(0)=z)}$.
The previous reasoning supplemented with this equality entails the desired statement.
Lemma \ref{L:main_equation} is proved.

The obtained equation has a natural interpretation. Namely, the term at the left-hand side of
\eqref{main_equation} accounts for the random walk trajectories starting at point $x$
and hitting point $y$ until time $t$. Every path of such kind belongs to one of two types.
The trajectories of the first type do not pass point $z$ before hitting $y$.
They are taken into account by the first summand at the right-hand side of \eqref{main_equation}.
As to a trajectory of the second type, it hits point $z$ at time $u$, $0\leq u\leq t$,
before reaching $y$ so that the part of the trajectory after hitting $z$ is a path starting
at $z$ and reaching $y$ until time $t-u$. The second summand at the right-hand side
of \eqref{main_equation} is responsible for the trajectories of the second type.

The next statement is a Lemma \ref{L:main_equation} corollary which provides the Laplace-Stieltjes
transform of a solution of equation \eqref{main_equation}.

\begin{Cr}\label{C:Laplace_transforms_of_Hxyz}
If $x,y,z\in\mathbb{Z}^{d}$ and $y\neq z$ then for any $\lambda>0$
\begin{equation}\label{LP_Hxyz}
\check{H}_{x,y,z}(\lambda)=\frac{\check{H}_{x,y}(\lambda)-\check{H}_{x,z}(\lambda)\check{H}_{z,y}(\lambda)}
{1-\check{H}_{z,y}(\lambda)\check{H}_{y,z}(\lambda)}.
\end{equation}
\end{Cr}
\textsc{Proof. }Due to Lemma \ref{L:main_equation} we have a system of two linear integral equations
in the functions $H_{x,y,z}(t)$ and $H_{x,z,y}(t)$
\begin{equation*}
\left\{
\begin{array}{lcl}
H_{x,y}(t)&=&H_{x,y,z}(t)+H_{x,z,y}\ast H_{z,y}(t),\\
H_{x,z}(t)&=&H_{x,z,y}(t)+H_{x,y,z}\ast H_{y,z}(t).
\end{array}
\right.
\end{equation*}
Application of the Laplace-Stieltjes transform to each equation of the system leads to a new system of algebraic
equations in $\check{H}_{x,y,z}(\lambda)$ and $\check{H}_{x,z,y}(\lambda)$
\begin{equation*}
\left\{
\begin{array}{lcl}
\check{H}_{x,y}(\lambda)&=&\check{H}_{x,y,z}(\lambda)+\check{H}_{x,z,y}(\lambda)\check{H}_{z,y}(\lambda),\\
\check{H}_{x,z}(\lambda)&=&\check{H}_{x,z,y}(\lambda)+\check{H}_{x,y,z}(\lambda)\check{H}_{y,z}(\lambda).
\end{array}
\right.
\end{equation*}
Solving this system we obtain \eqref{LP_Hxyz}. The proof is complete.

The next proposition clarifies the invariance of functions $H_{x,y,z}(\cdot)$ indexed by
$x,y,z\in\mathbb{Z}^{d}$, $z\neq y$, under shifts and reflections of space $\mathbb{Z}^{d}$.

\begin{Cr}\label{C:properties_of_Hxyz}
Let $r,x,y,z\in\mathbb{Z}^{d}$ and $y\neq z$. Then
\begin{equation}\label{properties_of_Hxyz}
H_{x+r,y+r,z+r}(t)=H_{x,y,z}(t),\quad H_{-x,-y,-z}(t)=H_{x,y,z}(t),\quad t\geq0.
\end{equation}
\end{Cr}
\textsc{Proof. }The desired assertion is a direct consequence of formulae \eqref{Hx,x=H0,0},
\eqref{Hx,y=Hx-y,0=Hy-x,0} and \eqref{LP_Hxyz}.

For proving Theorem 1 we also need the following lemma. Its derivation
is based on considering of all the possible jumps of the random walk after latency period $\tau$.

\begin{Lm}\label{L:H-xyz(infty)-H-xyz(t)=}
For $x,y,z\in\mathbb{Z}^{d}$ such that $y\neq z$, one has
\begin{equation}\label{H-xyz(infty)-H-xyz(t)=}
H^{-}_{x,y,z}(\infty)-H^{-}_{x,y,z}(t)=\sum\limits_{r\in\mathbb{Z}^{d},r\neq x, r\neq y, r\neq z}
{\frac{a(x,r)}{a}\left(H_{r,y,z}(\infty)-H_{r,y,z}(t)\right)},\quad t\geq0.
\end{equation}
\end{Lm}
\textsc{Proof. }In view of definitions of random variables $\tau$ and $\tau^{-}_{y,z}$ one has
\begin{eqnarray*}
H^{-}_{x,y,z}(\infty)-H^{-}_{x,y,z}(t)&=&{\sf P}(t<\tau^{-}_{y,z}<\infty|S(0)=x)\\
&=&\sum\limits_{r\in\mathbb{Z}^{d},r\notin\{x, y, z\}}{\!\!\!\!{\sf P}(t<\tau^{-}_{y,z}<\infty, S(\tau)=r|S(0)=x)}\\
&=&\sum\limits_{r\in\mathbb{Z}^{d},r\notin\{x, y, z\}}{\!\!\!\!{\sf P}(S(\tau)=r|S(0)=x)
{\sf P}(t<\tau^{-}_{y,z}<\infty|S(\tau)=r,S(0)=x)}.
\end{eqnarray*}
Here ${\sf P}(t<\tau^{-}_{y,z}<\infty|S(\tau)=r,S(0)=x)={\sf P}(t<\tau_{y,z}<\infty|S(0)=r)$
by the identity $\tau_{y,z}=\tau+\tau^{-}_{y,z}$ a.s. and due to strong Markov property of
the random walk w.r.t. the stopping time $\tau$ (see, e.g., \cite{Bremaud}, Ch.8, Sec.4).
Moreover, according to the explicit construction of $S$ the relation ${\sf P}(S(\tau)=r|S(0)=x)=a(x,r)/a$
is also true. Combining the obtained equalities with the previous argument
one comes to \eqref{H-xyz(infty)-H-xyz(t)=}. Lemma \ref{L:H-xyz(infty)-H-xyz(t)=} is proved.

The next result which is known as Fa\`{a} di Bruno's formula
plays an important role while establishing Theorem \ref{T:not_simple} for $d\geq3$.

\begin{Th}[see \cite{Faa_di_Bruno}, Appendix, Subsec.9]\label{Faa_di_Bruno}
Let functions $W(u)$ and $V(u)$ have $n$-th derivatives. Then for the
$n$-th derivative of function $U(u)=W(V(u))$ the following formula is valid
\begin{equation}\label{formula_Faa_di_Bruno}
U^{(n)}(u)=
\sum_{\substack{n_{1},n_{2},n_{3},\ldots\in\mathbb{Z}_{+}:\\
n_{1}+2n_{2}+3n_{3}+\ldots=n}}{W^{(n_{1}+n_{2}+n_{3}+\ldots)}
(V)\frac{n!}{n_{1}!\,n_{2}!\,n_{3}!\ldots}\left(\frac{V'}{1!}\right)^{n_{1}}
\left(\frac{V''}{2!}\right)^{n_{2}}\left(\frac{V'''}{3!}\right)^{n_{3}}\ldots}.
\end{equation}
\end{Th}

The last auxiliary result will be useful for proving Theorem \ref{T:simple}.

\begin{Lm}\label{L:trigonometric_equality}
For each $x\in\mathbb{N}$ the following equality holds true
\begin{equation}\label{trigonometric_equality}
\int\nolimits_{-\pi}^{\pi}{\frac{1-\cos{x\theta}}{1-\cos{\theta}}\,d\theta}=2\pi x.
\end{equation}
\end{Lm}
\textsc{Proof. }For each $n\in\mathbb{N}$ one has
$$\int\limits_{-\pi}^{\pi}{\frac{\cos{n\theta}-\cos{(n+1)\theta}}{1-\cos{\theta}}\,d\theta}\,
-\int\limits_{-\pi}^{\pi}{\frac{\cos{(n-1)\theta}-\cos{n\theta}}{1-\cos{\theta}}\,d\theta}$$
$$=\int\limits_{-\pi}^{\pi}{\frac{2\cos{n\theta}-\cos{(n+1)\theta}-\cos{(n-1)\theta}}{1
-\cos{\theta}}\,d\theta}=\int\limits_{-\pi}^{\pi}{\frac{2\cos{n\theta}(1-\cos{\theta})}{1
-\cos{\theta}}\,d\theta}=\int\limits_{-\pi}^{\pi}{2\cos{n\theta}\,d\theta}=0.$$
Whence for all $n\in\mathbb{N}$ the values of the following integrals are the same and equal,
for instance, the value of the integral for $n=1$
$$\int\nolimits_{-\pi}^{\pi}{\frac{\cos{(n-1)\theta}-\cos{n\theta}}{1-\cos{\theta}}\,d\theta}
=\int\nolimits_{-\pi}^{\pi}{\frac{1-\cos{\theta}}{1-\cos{\theta}}\,d\theta}=2\pi.$$
Hence, for each $x\in\mathbb{N}$ we obtain
$$\int\nolimits_{-\pi}^{\pi}{\frac{1-\cos{x\theta}}{1-\cos\theta}}=\sum\limits_{n=1}^{x}
{\int\nolimits_{-\pi}^{\pi}{\frac{\cos{(n-1)\theta}-\cos{n\theta}}{1-\cos{\theta}}\,d\theta}}
=\sum\limits_{n=1}^{x}{2\pi}=2\pi x.$$
Thus, equality \eqref{trigonometric_equality} is established and the proof of Lemma
\ref{L:trigonometric_equality} is complete.

\section{Proof of Theorem \ref{T:not_simple}}

Firstly, in view of Corollary \ref{C:properties_of_Hxyz} it suffices to prove Theorem
\ref{T:not_simple} for $z={\bf 0}$. Secondly, we have to consider separately the cases
$d=1$, $d=2$ and $d\geq3$, since the demonstration of the theorem essentially depends on
dimension $d$. Before passing to the non-simple random walk on $\mathbb{Z}$ let us write
down some formulae which have the same form for all $d\in\mathbb{N}$.

Let $x,y\in\mathbb{Z}^{d}\backslash\{{\bf 0}\}$ and $x\neq y$.
Using \eqref{universal_formula}, \eqref{lt_Hxx}, \eqref{lt_Hxy} and \eqref{LP_Hxyz}
as well as recalling that $\int\nolimits_{0}^{\infty}{e^{-\lambda t}(1-G(t))\, dt}=(\lambda+a)^{-1}$,
$\lambda\geq0$, we deduce the following relations
\begin{eqnarray}
\widehat{(H_{x,y,{\bf 0}}(\infty)-H_{x,y,{\bf
0}})}(\lambda)&=&\frac{H_{x,y,{\bf 0}}(\infty)}{\lambda}-\frac{G_{\lambda}(x,y)
G_{\lambda}({\bf 0},{\bf 0})-G_{\lambda}({\bf
0},x)G_{\lambda}({\bf 0},y)}{(G^{2}_{\lambda}({\bf 0},{\bf
0})-G^{2}_{\lambda}({\bf 0},y))\lambda},\label{lt_Hxy0(infty)-Hxy0}\\
\widehat{(H_{{\bf 0},y,{\bf 0}}(\infty)-H_{{\bf 0},y,{\bf
0}})}(\lambda)&=&\frac{H_{{\bf 0},y,{\bf 0}}(\infty)}{\lambda}-\frac{G_{\lambda}({\bf
0},y)}{(G^{2}_{\lambda}({\bf 0},{\bf
0})-G^{2}_{\lambda}({\bf 0},y))(\lambda+a)\lambda},\label{lt_Hoyo(infty)-Hoyo}\\
\widehat{(H_{y,y,{\bf 0}}(\infty)-H_{y,y,{\bf
0}})}(\lambda)&=&\frac{H_{y,y,{\bf 0}}(\infty)}{\lambda}-\frac{1}{\lambda}+\frac{G_{\lambda}({\bf
0},{\bf 0})}{(G^{2}_{\lambda}({\bf 0},{\bf
0})-G^{2}_{\lambda}({\bf 0},y))(\lambda+a)\lambda}.\label{lt_Hooy(infty)-Hooy}
\end{eqnarray}

\subsection{The case $d=1$}

In this subsection we focus on a non-simple random walk on $\mathbb{Z}$ although some arguments are
valid for a simple random walk on $\mathbb{Z}^{d}$ as well.
At the beginning we find the limit value $H_{x,y,{\bf 0}}(\infty)$ for $x,y\in\mathbb{Z}$, $y\neq{\bf 0}$,
due to the formula $H_{x,y,{\bf 0}}(\infty)=\lim\nolimits_{\lambda\to0+}{\check{H}_{x,y,{\bf 0}}(\lambda)}$.
To this end we write the following asymptotic representation
of the function $\check{H}_{{\bf 0},r}(\lambda)$, $r\in\mathbb{Z}$,
\begin{equation}\label{decomposition_d=1}
\check{H}_{{\bf 0},r}(\lambda)=1-\frac{\rho_{1}(r)
\,\sqrt{\lambda}}{a\,\gamma_{1}\sqrt{\pi}}+o(\sqrt{\lambda}),\quad\lambda\to0+,
\end{equation}
which results, in view of \eqref{universal_formula}, from application of Tauberian theorem
(Theorem 4 in \cite{Feller}, Ch.13, Sec.5) to relation \eqref{Hxy(infty)-Hxy(t)sim,d=1}.
Then substituting \eqref{decomposition_d=1} in \eqref{LP_Hxyz} and taking into account
\eqref{Hx,y=Hx-y,0=Hy-x,0} we calculate $\lim\nolimits_{\lambda\to0+}{\check{H}_{x,y,{\bf 0}}(\lambda)}$ and,
consequently, we get formula for $H_{x,y,{\bf 0}}(\infty)$ which is true for a simple random walk on
$\mathbb{Z}$ as well.
However, to complete the derivation of \eqref{T_Hxyz(infty),d=1,2}
for $d=1$ we have to show that $H_{x,y,{\bf 0}}(\infty)\in(0,1).$ Note that by virtue of
\eqref{Hxy(infty),d=1,d=2} equality \eqref{main_equation} implies the equivalences
$H_{x,y,{\bf 0}}(\infty)=0$ $\Leftrightarrow$ $H_{x,{\bf 0},y}(\infty)=1$ and
$H_{x,y,{\bf 0}}(\infty)=1$ $\Leftrightarrow$ $H_{x,{\bf 0},y}(\infty)=0$. Hence, to verify that
$H_{x,y,{\bf 0}}(\infty)\in(0,1)$ for all $x,y\in\mathbb{Z},$ $y\neq{\bf 0}$, it suffices to check that
$H_{x,y,{\bf 0}}(\infty)>0$ for all $x,y\in\mathbb{Z}$, $y\neq{\bf 0}$. In its turn the latter holds true if,
for instance, relation \eqref{H_xyz(infty)-H_xyz(t)sim,d=1} is satisfied with $C_{1}(x,y)>0$ for all
$x,y\in\mathbb{Z}$, $y\neq{\bf 0}$. Thus, the initial problem of validating that
$H_{x,y,{\bf 0}}(\infty)\in(0,1)$ is reduced to proving \eqref{H_xyz(infty)-H_xyz(t)sim,d=1}
and showing positivity of the function $C_{1}(\cdot,\cdot)$. The rest part of the subsection
is devoted to establishing the last two claims.

The demonstration of \eqref{H_xyz(infty)-H_xyz(t)sim,d=1} is based on formulae \eqref{lt_Hxy0(infty)-Hxy0},
\eqref{lt_Hoyo(infty)-Hoyo} and \eqref{lt_Hooy(infty)-Hooy} in which we substitute asymptotic decompositions
of functions $G_{\lambda}({\bf 0},{\bf 0})$ and $G_{\lambda}({\bf 0},{\bf 0})-G_{\lambda}({\bf 0},r)$, as
$\lambda\to0+$, for $r=x$, $r=y$ and $r=y-x$. Due to Tauberian theorem (Theorem 4 in \cite{Feller}, Ch.13, Sec.5)
the asymptotic behavior of $G_{\lambda}({\bf 0},{\bf 0})$ for $d=1$
is encoded in relation \eqref{p(t;x,y)sim}, namely,
\begin{equation}\label{G_lambda(0,0)sim_d=1}
G_{\lambda}({\bf 0},{\bf
0})=\frac{\gamma_{1}\sqrt{\pi}}{\sqrt{\lambda}}+o\left(\frac{1}{\sqrt{\lambda}}\right),\quad\lambda\to0+.
\end{equation}
Furthermore, using the same Tauberian theorem we infer that
\begin{equation}\label{G_lambda-G_lambda_asymptotic_decomposition_d=1}
G_{\lambda}({\bf 0},{\bf 0})-G_{\lambda}({\bf 0},r)=a^{-1}\rho_{1}(r)-2\sqrt{\pi}\,\tilde{\gamma}_{1}(r)\,
\sqrt{\lambda}+o(\sqrt{\lambda}),\quad\lambda\to0+,\quad r\in\mathbb{Z}\backslash\{{\bf 0}\}.
\end{equation}
Indeed, $a^{-1}\rho_{1}(r)=\int\nolimits_{0}^{\infty}{(p(t;{\bf
0},{\bf 0})-p(t;{\bf 0},r))\,dt}<\infty$ on account of \eqref{p(t;0,0)-p(t,x,y)sim} and
\begin{eqnarray}
G_{\lambda}({\bf 0},{\bf 0})-G_{\lambda}({\bf
0},r)&-&a^{-1}\rho_{1}(r)=\int\nolimits_{0}^{\infty}{e^{-\lambda
t}(p(t;{\bf 0},{\bf 0})-p(t;{\bf 0},r))\,dt}-a^{-1}\rho_{1}(r)\nonumber\\
&=&\lambda\int\nolimits_{0}^{\infty}{e^{-\lambda
t}\left(\int\nolimits_{0}^{t}{(p(u;{\bf 0},{\bf 0})-p(u;{\bf
0},r))\,du}\right)\,dt}-\lambda\int\nolimits_{0}^{\infty}{e^{-\lambda
t}a^{-1}\rho_{1}(r)\,dt}\nonumber\\
&=&-\lambda\int\nolimits_{0}^{\infty}{e^{-\lambda
t}\left(\int\nolimits_{t}^{\infty}{(p(u;{\bf 0},{\bf 0})-p(u;{\bf
0},r))\,du}\right)\,dt}\label{G_lambda-G_lambda_equality}.
\end{eqnarray}
Here $\int\nolimits_{t}^{\infty}{(p(u;{\bf 0},{\bf 0})-p(u;{\bf
0},r))\,du}\sim 2\tilde{\gamma}_{1}(r)\,t^{-1/2}$ (as $t\to\infty$) by formula \eqref{p(t;0,0)-p(t,x,y)sim}
and Theorem 31 in \cite{V_Lectures}. Now we substitute relation
\eqref{G_lambda-G_lambda_asymptotic_decomposition_d=1} and
the formula for $H_{x,y,{\bf 0}}(\infty)$ appearing in \eqref{T_Hxyz(infty),d=1,2}
into \eqref{lt_Hxy0(infty)-Hxy0}, \eqref{lt_Hoyo(infty)-Hoyo} and \eqref{lt_Hooy(infty)-Hooy}. After
collecting terms we substitute \eqref{G_lambda(0,0)sim_d=1} in the obtained formulae. Omitting tiresome
calculations we deduce that for $x,y\in\mathbb{Z},$ $y\neq{\bf 0}$, and $\lambda\to0+$
\begin{eqnarray*}
\widehat{(H_{x,y,{\bf 0}}(\infty)-H_{x,y,{\bf
0}})}(\lambda)&\sim&\frac{\rho_{1}(x)+\rho_{1}(y-x)-\rho_{1}(y)}{4\,a\,\sqrt{\pi}\,\gamma_{1}\sqrt{\lambda}}\\
&+&\frac{a\sqrt{\pi}\left(\rho_{1}(y)\tilde{\gamma}_{1}(x)+\rho_{1}(y-x)\tilde{\gamma}_{1}(y)-\rho_{1}(x)
\tilde{\gamma}_{1}(y)-\rho_{1}(y)\tilde{\gamma}_{1}(y-x)\right)}{\rho^{2}_{1}(y)\,\sqrt{\lambda}}.
\end{eqnarray*}
This can be rewritten in the form
\begin{equation}\label{lt_Hxy0(infty)-Hxy0_sim}
\widehat{(H_{x,y,{\bf 0}}(\infty)-H_{x,y,{\bf 0}})}(\lambda)
\sim\frac{C_{1}(x,y)\sqrt{\pi}}{\sqrt{\lambda}},\quad\lambda\to0+,
\end{equation}
given that $\tilde{\gamma}_{1}(r)=\pi r^{2} \gamma^{3}_{1},$ $r\in\mathbb{Z}$.
The latter identity is valid by virtue of definitions of $\gamma_{1}$ and  $\tilde{\gamma}_{1}(r)$ since
$\gamma_{1}=1/\sqrt{-2\pi\phi''(0)}$ and
$$\tilde{\gamma}_{1}(r)=\frac{r^{2}}{4\pi}\int\nolimits_{-\infty}^{+\infty}{\!\!\!\upsilon^{2}
e^{\phi''(0)\,\upsilon^{2}/2}\,d\upsilon}=\frac{r^{2}}{-4\pi\phi''(0)\sqrt{-\phi''(0)}}
\int\nolimits_{-\infty}^{+\infty}{\!\!\!u^{2}e^{-u^{2}/2}\,du}=\frac{r^{2}}{-2\phi''(0)\sqrt{-2\pi\phi''(0)}}.$$
We took into account that the variables change $u=\upsilon\sqrt{-\phi''(0)}$ reduces the integration
to writing $\int\nolimits_{-\infty}^{+\infty}{u^{2}e^{-u^{2}/2}\,du}=\sqrt{2\pi}$.
Note that relation \eqref{lt_Hxy0(infty)-Hxy0_sim} is valid for a simple random walk on $\mathbb{Z}$ as well.
However, in this case the function $C_{1}(x,y)$ takes nonnegative values (not strictly positive values).

If $C_{1}(x,y)>0$ then application of Tauberian theorem (Theorem 4 in \cite{Feller}, Ch.13, Sec.5)
to relation \eqref{lt_Hxy0(infty)-Hxy0_sim} leads to the desired statement \eqref{H_xyz(infty)-H_xyz(t)sim,d=1}.
Thus, to complete the proof of Theorem 1 for $d=1$ we only need to verify the positivity of function
$C_{1}(x,y)$, $x,y\in\mathbb{Z}$, $y\neq{\bf 0}$.

According to Corollary \ref{C:properties_of_Hxyz} we may assume that $y>{\bf 0}$.
Then there exist three possible relative positions of points $x,y$ and ${\bf 0}$, namely,
$x\geq y$, ${\bf 0}\leq x<y$
and $x<{\bf 0}$. Let us consider at first the case $x\geq y$. To check strict positivity of $C_{1}(x,y)$
for such $x$ and $y$ we search for a lower estimate of it. To this end note that
$\rho_{1}(x)+\rho_{1}(y-x)-\rho_{1}(y)\geq0$ in view of formula for $H_{y,x,{\bf 0}}(\infty)$ appearing
in \eqref{T_Hxyz(infty),d=1,2} and the evident inequality $H_{y,x,{\bf 0}}(\infty)\leq 1$.
Obviously, $C_{1}({\bf 0},y)\geq0$ and, consequently,
$\rho^{2}_{1}(y)-4 a^{2} \pi^{2} y^{2} \gamma^{4}_{1}\geq0$.
Combining these inequalities we see that for $x\geq y$
\begin{eqnarray*}
C_{1}(x,y)&=&\frac{(\rho_{1}(x)+\rho_{1}(y-x)-\rho_{1}(y))(\rho^{2}_{1}(y)-4 a^{2} \pi^{2} y^{2} \gamma_{1}^{4})}
{4\,a\,\pi\,\rho^{2}_{1}(y)\,\gamma_{1}}+\frac{2\rho_{1}(y-x)a \pi y^{2} \gamma^{3}_{1}}{\rho^{2}_{1}(y)}\\
&+&\frac{2 a \pi \gamma^{3}_{1} y(x-y)}{\rho_{1}(y)}\geq\frac{2\rho_{1}(y-x)a \pi y^{2} \gamma^{3}_{1}}{\rho^{2}_{1}(y)}
+\frac{2 a \pi \gamma^{3}_{1} y(x-y)}{\rho_{1}(y)}>0.
\end{eqnarray*}

Before passing to the case ${\bf 0}\leq x<y$ we derive one more useful relation. Recall that
$\tau_{y,{\bf 0}}\geq\tau^{-}_{y,{\bf 0}}$ a.s. because $\tau_{y,{\bf 0}}=\tau^{-}_{y,{\bf 0}}+\tau$ a.s.
Since the random variable $\tau$ is exponentially distributed,
${\sf P}(\tau_{y,{\bf 0}}<\infty)={\sf P}(\tau^{-}_{y,{\bf 0}}<\infty)$.
Therefore ${\sf P}(t<\tau_{y,{\bf 0}}<\infty)\geq{\sf P}(t<\tau^{-}_{y,{\bf 0}}<\infty)$ and hence
$H_{x,y,{\bf 0}}(\infty)-H_{x,y,{\bf 0}}(t)\geq H^{-}_{x,y,{\bf 0}}(\infty)-H^{-}_{x,y,{\bf 0}}(t)$.
Due to Lemma \ref{L:H-xyz(infty)-H-xyz(t)=} and the monotonicity property of the Laplace transform
it follows that
\begin{equation}\label{Lp_H-xyz(infty)-H-xyz(t)=}
\widehat{(H_{x,y,{\bf 0}}(\infty)-H_{x,y,{\bf 0}})}(\lambda)\geq\sum\limits_{r\in\mathbb{Z}, r\neq x, r\neq y,
r\neq{\bf 0}}{\frac{a(x,r)}{a}\widehat{(H_{r,y,{\bf 0}}(\infty)-H_{r,y,{\bf 0}})}(\lambda)},\quad \lambda>0.
\end{equation}

Now turn to the case ${\bf 0}\leq x<y$. We show the positivity of $C_{1}(x,y)$ successively for
$x=y-1$, $x=y-2$, $\ldots$, $x={\bf 0}$. If $x=y-1$ then there exists a point $r>y$ such that
$a(y-1,r)>0$ as otherwise the random walk is simple.
Therefore, by virtue of \eqref{lt_Hxy0(infty)-Hxy0_sim} and of the just established inequality $C_{1}(r,y)>0$, ${\bf 0}<y<r$,
there is a summand of the order $1/\sqrt{\lambda}$ at the right-hand side of \eqref{Lp_H-xyz(infty)-H-xyz(t)=}.
Consequently, by \eqref{lt_Hxy0(infty)-Hxy0_sim} the left-hand side of \eqref{Lp_H-xyz(infty)-H-xyz(t)=}
is also of the order $1/\sqrt{\lambda}$, $\lambda\to0+$, and the relation $C_{1}(y-1,y)>0$ with $y>{\bf 0}$
is proved. Recall that the random walk under
consideration is irreducible, that is, the greatest common divisor of all $r\in\mathbb{Z}\backslash\{{\bf 0}\}$
such that $a(r)>0$ equals 1 (see \cite{Borovkov}, Ch.12, Sec.3). This implies for $x=y-2$ that there exists
a point $r>y$ or $r=y-1$
such that $a(y-2,r)>0$. Since for such $r$ we have just shown the positivity of $C_{1}(r,y)$, in view of
\eqref{lt_Hxy0(infty)-Hxy0_sim} it means that the left-hand side of \eqref{Lp_H-xyz(infty)-H-xyz(t)=} has the order
of growth $1/\sqrt{\lambda}$, as $\lambda\to0$. So we conclude that
$C_{1}(y-2,y)>0$ with $y>{\bf 0}$. Considering successively the cases $x=y-3$, $x=y-4$, $\ldots$, $x={\bf 0}$,
we prove according to the described scheme that $C_{1}(x,y)>0$ for all $x\in[{\bf 0},y),$ $x\in\mathbb{Z}$.

The case $x<{\bf 0}<y$ is also treated by indirect methods, since the immediate estimation of $C_{1}(x,y)$ for
such $x$ and $y$ is rather difficult. We consider in succession $x=-1$, $x=-2$, $\ldots$. At first we assume
that for the random walk $S$ there exist at least three
points $r_{1}>{\bf 0}$, $r_{2}>{\bf 0}$ and $r_{3}>{\bf 0}$ such that $a(r_{1})>0$, $a(r_{2})>0$ and
$a(r_{3})>0$. If $x=-1$ then at the right-hand side of \eqref{Lp_H-xyz(infty)-H-xyz(t)=}
there is at least one summand indexed by $r$ where $r>{\bf 0}$, $r\neq y$ and $a(-1,r)>0$.
Similarly to the previous discussion we infer that the left-hand side of \eqref{Lp_H-xyz(infty)-H-xyz(t)=}
has the order of growth $1/\sqrt{\lambda}$, as $\lambda\to0+$, and $C_{1}(-1,y)>0$.
If $x=-2$ then at the right-hand side of \eqref{Lp_H-xyz(infty)-H-xyz(t)=}
there exists at least one summand indexed by $r=-1$ or $r>{\bf 0}$, such that $r\neq y$ and $a(-2,r)>0$.
Since for such $r$ the positivity of $C_{1}(r,y)$ has been already proved, in view of
\eqref{lt_Hxy0(infty)-Hxy0_sim} and \eqref{Lp_H-xyz(infty)-H-xyz(t)=} we deduce that $C_{1}(-2,y)>0$.
The positivity of $C_{1}(x,y)$ for $x=-3$, $x=-4$, $\ldots$ is verified in the same manner.

We continue to deal with the case $x<{\bf 0}<y$ assuming that there exist exactly two points
$r_{1}>{\bf 0}$ and $r_{2}>{\bf 0}$ such that $a(r_{1})>0$ and $a(r_{2})>0$.
Let $x=-1$. If there is a point $r>{\bf 0}$, $r\neq y$, such that $a(-1,r)>0$ (i.e.
$r_{1}=r+1$ or $r_{2}=r+1$) then by virtue of \eqref{lt_Hxy0(infty)-Hxy0_sim} and
\eqref{Lp_H-xyz(infty)-H-xyz(t)=} as well as of the verified inequality $C_{1}(r,y)>0$ we get the desired
relation $C_{1}(-1,y)>0$. Otherwise, one has $r_{1}=1$ and $r_{2}=y+1$. If $y>1$ then $a(-2,y-1)=a(y+1)>0$ and,
consequently, $C_{1}(-2,y)>0$ in view of \eqref{lt_Hxy0(infty)-Hxy0_sim}, \eqref{Lp_H-xyz(infty)-H-xyz(t)=} and
the established earlier estimate $C_{1}(y-1,y)>0$. Therefore, if $r_{1}=1$, $r_{2}=y+1$ and $y>1$
we obtain $C_{1}(-1,y)>0$ due to \eqref{lt_Hxy0(infty)-Hxy0_sim} and \eqref{Lp_H-xyz(infty)-H-xyz(t)=}
as well as the relation $a(-1,-2)=a(1)>0$ and the just proved inequality $C_{1}(-2,y)>0$.
However, for $r_{1}=1$, $r_{2}=y+1$ and $y=1$, the previous arguments do not work and
we have to employ the following facts concerning the embedded chain $\{S_{n}, n\in\mathbb{Z}_{+}\}$.

For the embedded chain $\{S_{n}, n\in\mathbb{Z}_{+}\}$ introduce a passage time $T:=\min\{n>0: S_{n}\geq -1\}$.
According to \cite{Afanasiev}, Theorems 15.1 and 15.2,
\begin{equation}\label{P(T>n)sim}
{\sf P}(T<\infty|S(0)=-1)=1,\quad{\sf P}(T>n|S(0)=-1)\sim\frac{c}{\sqrt{n}},\quad n\to\infty,
\end{equation}
where $c$ is some positive constant.

For brevity we write ${\sf P}_{-1}(\cdot):={\sf P}(\cdot|S(0)=-1)$.
Since we assume that $r_{1}=1$ and $r_{2}=2$, the process $\{S_{n}, n\in\mathbb{Z}_{+}\}$
can perform jumps to the nearest-neighbor points or to the next nearest points. Hence,
\begin{eqnarray}\label{P(T=n,S_n=-1)>=}
{\sf P}_{-1}(T=n, S_{n}=-1)\!\!&=&\!\!{\sf P}_{-1}(T=n, S_{n-1}=-2, S_{n}=-1)+{\sf P}_{-1}(T=n, S_{n-1}=-3,S_{n}=-1)\nonumber\\
\!\!&=&\!\!{\sf P}_{-1}(S_{k}<-1, 0<k<n-1, S_{n-1}=-2, S_{n}=-1)\nonumber\\
\!\!&+&\!\!{\sf P}_{-1}(S_{k}<-1, 0<k<n-1, S_{n-1}=-3, S_{n}=-1)\nonumber\\
\!\!&=&\!\!\frac{a(-2,-1)}{a}{\sf P}_{-1}(S_{k}<-1, 0<k<n-1, S_{n-1}=-2)\nonumber\\
\!\!&+&\!\!\frac{a(-3,-1)}{a}{\sf P}_{-1}(S_{k}<-1, 0<k<n-1, S_{n-1}=-3)\nonumber\\
\!\!&\geq&\!\!\frac{\min\{a(1),a(2)\}}{a}{\sf P}_{-1}(T=n).
\end{eqnarray}

Let us return to a lower estimate of $H_{-1,1,{\bf 0}}(\infty)-
H_{-1,1,{\bf 0}}(t)={\sf P}_{-1}(t<\tau_{1,{\bf 0}}<\infty)$, as $t\to\infty$.
Denote by $N=\{N(t), t\geq0\}$ the Poisson process constructed by means of the random sequence
$\{\tau^{(n+1)}-\tau^{(n)}\}_{n=0}^{\infty}$, i.e. $N$ is the Poisson process with intensity $a$.
Considering all the possible jumps of $S$, taking into account \eqref{P(T=n,S_n=-1)>=}
and invoking the explicit construction of $S$ we derive that
\begin{eqnarray}\label{P(t<tau_10<infty)>=}
{\sf P}_{-1}(t<\tau_{1,{\bf 0}}<\infty)&=&\sum\limits_{n=1}^{\infty}{{\sf P}_{-1}(\tau_{1,{\bf 0}}=\tau^{(n)},
t<\tau_{1,{\bf 0}}<\infty)}\nonumber\\
&=&\sum\limits_{n=1}^{\infty}{{\sf P}_{-1}(S_{k}\neq{\bf 0}, S_{k}\neq1, 0<k<n, S_{n}=1,
t<\tau^{(n)}<\infty)}\nonumber\\
&=&\sum\limits_{n=1}^{\infty}{{\sf P}_{-1}(S_{k}\leq-1, 0<k<n-1, S_{n-1}=-1, S_{n}=1, N(t)<n)}\nonumber\\
&=&\sum\limits_{n=1}^{\infty}{{\sf P}_{-1}(S_{k}\leq-1, 0<k<n-1, S_{n-1}=-1, S_{n}=1) {\sf P}(N(t)<n)}\nonumber\\
&\geq&\frac{a(-1,1)}{a}\sum\limits_{n=1}^{\infty}{{\sf P}_{-1}(S_{k}<-1, 0<k<n-1, S_{n-1}=-1)
{\sf P}(N(t)<n)}\nonumber\\
&=&\frac{a(2)}{a}\sum\limits_{n=1}^{\infty}{{\sf P}_{-1}(T=n-1, S_{n-1}=-1){\sf P}(N(t)<n)}\nonumber\\
&\geq&\frac{a(2)\min\{a(1),a(2)\}}{a^{2}}\sum\limits_{n=1}^{\infty}{{\sf P}_{-1}(T=n-1)\sum\limits_{k=0}^{n-1}{\frac{e^{-at}(at)^{k}}{k!}}}\nonumber\\
&=&\frac{a(2)\min\{a(1),a(2)\}}{a^{2}}\sum\limits_{k=0}^{\infty}{{\sf P}_{-1}(T\geq k)\frac{e^{-at}(at)^{k}}{k!}}.
\end{eqnarray}
Thus, finding a lower estimate of ${\sf P}_{-1}(t<\tau_{1,{\bf 0}}<\infty)$ is reduced to
establishing a lower estimate of the last series.
It is not difficult to check that $(N(t)-at)/\sqrt{at}\stackrel{Law}\longrightarrow\xi$ with $\xi\sim\mathcal{N}(0,1)$,
as $t\to\infty$. Consequently,
${\sf P}(N(t)\in(a t-b\sqrt{at},a t+b\sqrt{a t}))\to
{\sf P}(\xi\in(-b,b))>0$, as $t\to\infty$, for each $b>0$. Hence, in view of \eqref{P(T>n)sim}
\begin{eqnarray}\label{series_sim}
\sum\limits_{k=0}^{\infty}{{\sf P}_{-1}(T\geq k)\frac{e^{-at}(at)^{k}}{k!}}&\geq&{\sf P}(T\geq[at+b\sqrt{at}])\,
{\sf P}(N(t)\in(at-b\sqrt{at},at+b\sqrt{at}))\nonumber\\
&\sim&\frac{c}{\sqrt{at}}\,{\sf P}(\xi\in(-b,b)),\quad t\to\infty.
\end{eqnarray}
Combining \eqref{lt_Hxy0(infty)-Hxy0_sim}, \eqref{P(t<tau_10<infty)>=} and \eqref{series_sim} as well as
taking into account the monotonicity property of the Laplace transform we conclude that $C_{1}(-1,1)>0$.

Thus, we have shown that $C_{1}(-1,y)>0$ for each $y>0$. Let us pass to the starting point $x=-2$.
Recall that we assume that there exist exactly two points $r_{1}>{\bf 0}$ and $r_{2}>{\bf 0}$
such that $a(r_{1})>0$ and $a(r_{2})>0$. If there is point $r>{\bf 0}$, $r\neq y$, or $r=-1$ such that
$a(-2,r)>0$ (that is, $r_{1}=r+2$ or $r_{2}=r+2$) then due to \eqref{lt_Hxy0(infty)-Hxy0_sim} and
\eqref{Lp_H-xyz(infty)-H-xyz(t)=} as well as the verified inequality $C_{1}(r,y)>0$ we get the desired
estimate $C_{1}(-2,y)>0$. Otherwise, one has $r_{1}=2$ and $r_{2}=y+2$. In particular, $a(-4,-2)>0$ and
$a(-4,y-2)>0$. Since the random walk is irreducible, one can guarantee that $y-2\neq{\bf 0}$.
Consequently, $C_{1}(-4,y)>0$ in view of \eqref{lt_Hxy0(infty)-Hxy0_sim}, \eqref{Lp_H-xyz(infty)-H-xyz(t)=}
and the earlier established inequality $C_{1}(y-2,y)>0$ with $y>{\bf 0}$ and $y\neq 2$.
In its turn $C_{1}(-2,y)>0$ by virtue of \eqref{lt_Hxy0(infty)-Hxy0_sim}, \eqref{Lp_H-xyz(infty)-H-xyz(t)=}
and the just proved estimate $C_{1}(-4,y)>0$. Therefore, the positivity of $C_{1}(x,y)$ for $x=-2$ and $y>{\bf 0}$
is demonstrated. For $x=-3$, $x=-4$, $\ldots$, verification of $C_{1}(x,y)$ positivity, $y>{\bf 0}$,
is implemented similarly to the case $x=-2$. Consequently, Theorem~\ref{T:not_simple} for $d=1$ is proved.

\subsection{The case $d=2$}

As in the previous subsection one can find the limit value $H_{x,y,{\bf 0}}(\infty)$ for $x,y\in\mathbb{Z}^{2}$,
$y\neq{\bf 0}$, by the formula
$H_{x,y,{\bf 0}}(\infty)=\lim\nolimits_{\lambda\to0+}{\check{H}_{x,y,{\bf 0}}(\lambda)}$.
For this purpose we write the following asymptotic representation
of the function $\check{H}_{{\bf 0},r}(\lambda)$, $r\in\mathbb{Z}^{2}$,
\begin{equation}\label{decomposition_d=2}
\check{H}_{{\bf 0},r}(\lambda)=1+\frac{\rho_{2}(r)}
{a\,\gamma_{2}\ln{\lambda}}+o\left(\frac{1}{\ln{\lambda}}\right),\quad\lambda\to0+,
\end{equation}
which results from application of Tauberian theorem (Theorem 4 in \cite{Feller}, Ch.13, Sec.5) to relation
\eqref{Hxy(infty)-Hxy(t)sim,d=2} and from employing \eqref{universal_formula}.
Substituting \eqref{decomposition_d=2} in \eqref{LP_Hxyz} and taking into account
\eqref{Hx,y=Hx-y,0=Hy-x,0} we find $\lim\nolimits_{\lambda\to0+}{\check{H}_{x,y,{\bf 0}}(\lambda)}$ and,
thus, $H_{x,y,{\bf 0}}(\infty)$ is obtained. However, to complete the derivation of \eqref{T_Hxyz(infty),d=1,2}
for $d=2$ one has to show that $H_{x,y,{\bf 0}}(\infty)\in(0,1).$ Similarly to Subsection 4.1,
equalities \eqref{Hxy(infty),d=1,d=2} and \eqref{main_equation} entail that
$H_{x,y,{\bf 0}}(\infty)=0$ $\Leftrightarrow$ $H_{x,{\bf 0},y}(\infty)=1$ and
$H_{x,y,{\bf 0}}(\infty)=1$ $\Leftrightarrow$ $H_{x,{\bf 0},y}(\infty)=0$. Therefore, for
establishing that $H_{x,y,{\bf 0}}(\infty)\in(0,1)$ for all $x,y\in\mathbb{Z}^{2},$ $y\neq{\bf 0}$,
it suffices to verify the inequality $H_{x,y,{\bf 0}}(\infty)>0$ for all $x,y\in\mathbb{Z}^{2}$,
$y\neq{\bf 0}$. In its turn the latter holds true if,
for instance, relation \eqref{H_xyz(infty)-H_xyz(t)sim,d=2} is satisfied with $C_{2}(x,y)>0$ for all
$x,y\in\mathbb{Z}^{2}$, $y\neq{\bf 0}$. Hence, the initial problem of checking that
$H_{x,y,{\bf 0}}(\infty)\in(0,1)$ is reduced to proving of \eqref{H_xyz(infty)-H_xyz(t)sim,d=2}
and showing the positivity of function $C_{2}(\cdot,\cdot)$. The rest part of this subsection
is devoted to validating the last two claims.

The proof of \eqref{H_xyz(infty)-H_xyz(t)sim,d=2} bears on formulae
\eqref{lt_Hxy0(infty)-Hxy0}--\eqref{lt_Hooy(infty)-Hooy} in which we substitute asymptotic decompositions
of functions $G_{\lambda}({\bf 0},{\bf 0})$ and $G_{\lambda}({\bf 0},{\bf 0})-G_{\lambda}({\bf 0},r)$, as
$\lambda\to0+$, for $r=x$, $r=y$ and $r=y-x$. Due to Tauberian theorem (Theorem 2 in \cite{Feller}, Ch.13, Sec.5)
the asymptotic behavior of $G_{\lambda}({\bf 0},{\bf 0})$ for $d=2$ is implied by \eqref{p(t;x,y)sim},
namely,
\begin{equation}\label{G_lambda(0,0)sim_d=2}
G_{\lambda}({\bf 0},{\bf
0})=-\gamma_{2}\ln{\lambda}+o(\ln{\lambda}),\quad\lambda\to0+.
\end{equation}
Using the same Tauberian theorem and relation \eqref{G_lambda-G_lambda_equality}, valid for
all $d\in\mathbb{N}$, we deduce that
\begin{equation}\label{G_lambda-G_lambda_asymptotic_decomposition_d=2}
G_{\lambda}({\bf 0},{\bf 0})-G_{\lambda}({\bf 0},r)=a^{-1}\rho_{2}(r)+\tilde{\gamma}_{2}(r)\,
\lambda\ln{\lambda}+o(\lambda\ln{\lambda}),\quad\lambda\to0+,\quad r\in\mathbb{Z}^{2}\backslash\{{\bf 0}\}.
\end{equation}
Indeed, in view of \eqref{p(t;0,0)-p(t,x,y)sim} one has $a^{-1}\rho_{2}(r)=\int\nolimits_{0}^{\infty}{(p(t;{\bf
0},{\bf 0})-p(t;{\bf 0},r))\,dt}<\infty$ and by formula \eqref{p(t;0,0)-p(t,x,y)sim}
and Theorem 31 in \cite{V_Lectures} one infers that $\int\nolimits_{t}^{\infty}{(p(u;{\bf 0},{\bf 0})-p(u;{\bf
0},r))\,du}\sim\tilde{\gamma}_{2}(r)/t$, as $t\to\infty$. Now we substitute relation
\eqref{G_lambda-G_lambda_asymptotic_decomposition_d=2} and
the formula for $H_{x,y,{\bf 0}}(\infty)$ appearing in \eqref{T_Hxyz(infty),d=1,2}
into \eqref{lt_Hxy0(infty)-Hxy0}--\eqref{lt_Hooy(infty)-Hooy}. After
collecting terms we substitute \eqref{G_lambda(0,0)sim_d=2} in the obtained formulae. Omitting tedious
calculations we establish that
$$\widehat{(H_{x,y,{\bf 0}}(\infty)-H_{x,y,{\bf 0}})}(\lambda)\sim\frac{\rho_{2}(x)+\rho_{2}(y-x)-
\rho_{2}(y)}{-4\,a\,\gamma_{2}\,\lambda\ln{\lambda}},\quad\lambda\to0+,\quad x,y\in\mathbb{Z}^{2},
\quad y\neq{\bf 0}.$$
If $C_{2}(x,y)>0$ then application of Tauberian theorem (Theorem 4 in \cite{Feller}, Ch.13, Sec.5)
to the last asymptotic relation gives \eqref{H_xyz(infty)-H_xyz(t)sim,d=2}.
Thus, to complete the proof of Theorem \ref{T:not_simple} for $d=2$ we only have to verify the positivity
of $C_{2}(x,y)$ for all $x,y\in\mathbb{Z}^{2}$, $y\neq{\bf 0}$.
To prove the claim we employ the method of ad absurdum. Assume that there exist
$x$ and $y$ such that $x\in\mathbb{Z}^{2}$, $y\in\mathbb{Z}^{2}\backslash\{{\bf 0}\}$ and
$C_{2}(x,y)=0$. Then the definition of $C_{2}(x,y)$ implies that
\begin{equation}\label{C2(x,y)=0}
\rho_{2}(x)+\rho_{2}(y-x)-\rho_{2}(y)=0.
\end{equation}
In view of \eqref{T_Hxyz(infty),d=1,2} the last equality is equivalent to the
relation $H_{y,x,{\bf 0}}(\infty)=1$. Therefore, formulae \eqref{Hxy(infty),d=1,d=2} and \eqref{main_equation}
entail that $H_{y,{\bf 0},x}(t)=0$ and $H_{y,x,{\bf 0}}(t)=H_{y,x}(t)$ for $t\geq0$.
Hence, by virtue of \eqref{Hxy(infty)-Hxy(t)sim,d=2} one has $C_{2}(y,x)>0$ and
$C_{2}(y,x)=\rho_{2}(y-x)/(a\,\gamma_{2})$. Taking into account the definition of $C_{2}(y,x)$
the last equality can be rewritten as $\rho_{2}(y)-\rho_{2}(x)=3\rho_{2}(y-x)$. This relation
combined with assumption \eqref{C2(x,y)=0} leads to contradiction. Hence, $C_{2}(x,y)>0$ for all
$x,y\in\mathbb{Z}^{2}$, $y\neq{\bf 0}$. The proof of Theorem \ref{T:not_simple} for $d=2$ is complete.

\subsection{The case $d\geq3$}

In this subsection we assume that $d\geq3$.
Basing on the formula $H_{x,y,{\bf 0}}(\infty)=\lim\nolimits_{\lambda\to0+}{\check{H}_{x,y,{\bf 0}}(\lambda)}$,
$x,y\in\mathbb{Z}^{d},$ $y\neq{\bf 0}$,
we find the limit value $H_{x,y,{\bf 0}}(\infty)$ with the help of
\eqref{Hxx(infty),d>=3}, \eqref{Hxy(infty),d>=3} and \eqref{LP_Hxyz}. To prove \eqref{T_Hxyz(infty),d>=3}
completely we have to check that $H_{x,y,{\bf 0}}(\infty)\in(0,1)$. Obviously, in view of
\eqref{Hxx(infty),d>=3}, \eqref{Hxy(infty),d>=3} and \eqref{main_equation} the relation
$H_{x,y,{\bf 0}}(\infty)\leq H_{x,y}(\infty)<1$
holds true. For verifying inequality $H_{x,y,{\bf 0}}(\infty)>0$ it suffices to establish
\eqref{H_xyz(infty)-H_xyz(t)sim,d>=3} and show the positivity of $C_{d}(x,y)$.
So the rest part of the subsection is devoted to validating these two assertions.

Further on we suppose that $x$ and $y$
are fixed points within $\mathbb{Z}^{d}\backslash\{{\bf 0}\}$.
Formulae \eqref{lt_Hxy0(infty)-Hxy0} -- \eqref{lt_Hooy(infty)-Hooy} can be rewritten in the following form
$$\widehat{(H_{x,y,{\bf 0}}(\infty)-H_{x,y,{\bf 0}})}(\lambda)
=\frac{\rho_{d}(y)(\rho_{d}(y-x)+\rho_{d}(x)-\rho_{d}(y))}
{a^{2}(G_{0}({\bf 0},{\bf 0})+G_{0}({\bf 0},y))}J_{1}(\lambda;{\bf
0},y)+J_{2}(\lambda;x,y)+J_{3}(\lambda;x,y)$$
\begin{equation}\label{lt_Hxy0(infty)-Hxy0_d>=3_sum}
-\frac{G_{0}({\bf 0},{\bf 0})(\rho_{d}(y-x)+\rho_{d}(x)-\rho_{d}(y))}{a(G_{0}({\bf 0},{\bf
0})+G_{0}({\bf 0},y))}J_{4}(\lambda;y)-\frac{a(G_{0}({\bf 0},{\bf 0})G_{0}(x,y)-G_{0}({\bf 0},x)G_{0}({\bf
0},y))}{\rho_{d}(y)(G_{0}({\bf 0},{\bf 0})+G_{0}({\bf
0},y))}J_{5}(\lambda;y),
\end{equation}
\begin{eqnarray}
\widehat{(H_{{\bf 0},y,{\bf 0}}(\infty)-H_{{\bf 0},y,{\bf
0}})}(\lambda)&=&\frac{\rho_{d}(y)}{a^{2}(G_{0}({\bf 0},{\bf 0})+G_{0}({\bf
0},y))}J_{1}(\lambda;y;y)-\frac{G_{0}({\bf 0},y)}{a(G_{0}({\bf
0},{\bf 0})+G_{0}({\bf 0},y))}J_{4}(\lambda;y)\nonumber\\
&-&\frac{G_{0}({\bf 0},y)}{\rho_{d}(y)(G_{0}({\bf 0},{\bf
0})+G_{0}({\bf
0},y))}J_{5}(\lambda;y)+\frac{1}{a}J_{6}(\lambda;y,y),\label{lt_Hoyo(infty)-Hoyo_d>=3_sum}
\end{eqnarray}
\begin{eqnarray}
\widehat{(H_{y,y,{\bf 0}}(\infty)-H_{y,y,{\bf 0}})}(\lambda)&
=&\frac{\rho_{d}(y)}{a^{2}(G_{0}({\bf 0},{\bf 0})+G_{0}({\bf 0},y))}J_{1}(\lambda;{\bf
0},y)-\frac{G_{0}({\bf 0},{\bf 0})}{a(G_{0}({\bf 0},{\bf
0})+G_{0}({\bf 0},y))}J_{4}(\lambda;y)\nonumber\\
&+&\frac{G_{0}({\bf 0},{\bf 0})}{\rho_{d}(y)(G_{0}({\bf 0},{\bf
0})+G_{0}({\bf
0},y))}J_{5}(\lambda;y)-\frac{1}{a}J_{6}(\lambda;{\bf
0},y)\label{lt_Hooy(infty)-Hooy_d>=3_sum}
\end{eqnarray}
where for $r\in\mathbb{Z}^{d}$ and $\lambda>0$
\begin{eqnarray*}
J_{1}(\lambda;r,y):=\frac{K_{d}(\lambda;r)}{G^{2}_{\lambda}({\bf
0},{\bf 0})-G^{2}_{\lambda}({\bf
0},y)},&\;&J_{2}(\lambda;r,y):=\frac{(K_{d}(\lambda;y-r)-K_{d}(\lambda;y))G_{\lambda}
({\bf 0},{\bf 0})}{G^{2}_{\lambda}({\bf 0},{\bf 0})-G^{2}_{\lambda}({\bf 0},y)},\\
J_{3}(\lambda;r,y):=\frac{(K_{d}(\lambda;{\bf 0})-K_{d}(\lambda;r))G_{\lambda}({\bf 0},y)}
{G^{2}_{\lambda}({\bf 0},{\bf 0})-G^{2}_{\lambda}({\bf 0},y)},&\;&
J_{4}(\lambda;y):=\frac{K_{d}(\lambda;{\bf
0})-K_{d}(\lambda;y)}{G^{2}_{\lambda}({\bf 0},{\bf
0})-G^{2}_{\lambda}({\bf 0},y)},\\
J_{5}(\lambda;y):=\frac{K_{d}(\lambda;{\bf
0})-K_{d}(\lambda;y)}{G_{\lambda}({\bf 0},{\bf 0})-G_{\lambda}({\bf
0},y)},&\;&
J_{6}(\lambda;r,y):=\frac{G_{\lambda}({\bf
0},r)}{(\lambda+a)(G^{2}_{\lambda}({\bf 0},{\bf
0})-G^{2}_{\lambda}({\bf 0},y))}
\end{eqnarray*}
and
$$K_{d}(\lambda;r):=\frac{G_{0}({\bf 0},r)-G_{\lambda}({\bf
0},r)}{\lambda}.$$

Let us briefly describe the further scheme of the proof of \eqref{H_xyz(infty)-H_xyz(t)sim,d>=3}.
At first we differentiate $[d/2]-1$ times (symbol $[\cdot]$ stands for the integer part of a number)
the both sides of \eqref{lt_Hxy0(infty)-Hxy0_d>=3_sum},
\eqref{lt_Hoyo(infty)-Hoyo_d>=3_sum} and \eqref{lt_Hooy(infty)-Hooy_d>=3_sum} w.r.t.
variable $\lambda$ at point $\lambda>0$. Then we retrieve the asymptotic behavior of function
$J^{([d/2]-1)}_{1}(\lambda;r,y)$, as $\lambda\to0+$, for each ${r\in\mathbb{Z}^{d}}$.
Moreover, we verify that the functions $J^{([d/2]-1)}_{2}(\lambda;r,y),$ $J^{([d/2]-1)}_{3}(\lambda;r,y),$
$J^{([d/2]-1)}_{4}(\lambda;y),$ $J^{([d/2-1])}_{5}(\lambda;y)$ and $J^{([d/2]-1)}_{6}(\lambda;r,y)$
are $o(J^{([d/2]-1)}_{1}(\lambda;r,y))$, as $\lambda\to0+$, for each $r$.
Afterwards one can employ Corollary 43 in \cite{V_Lectures} and, thus, find the asymptotic behavior
of the functions $H_{x,y,{\bf 0}}(\infty)-H_{x,y,{\bf 0}}(t)$,
${H_{{\bf 0},y,{\bf 0}}(\infty)-H_{{\bf 0},y,{\bf 0}}(t)}$ and $H_{y,y,{\bf 0}}(\infty)-H_{y,y,{\bf 0}}(t)$,
as $t\to\infty$.

Before studying the limiting behavior of $J^{([d/2]-1)}_{1}(\lambda;r,y)$
we have to investigate the asymptotic properties of the functions
$K_{d}(\lambda;r)$, $(G_{\lambda}({\bf 0},{\bf 0})-G_{\lambda}({\bf
0},y))^{-1}$, $(G^{2}_{\lambda}({\bf 0},{\bf
0})-G^{2}_{\lambda}({\bf 0},y))^{-1}$, $G_{\lambda}({\bf 0},r)/(G^{2}_{\lambda}({\bf 0},{\bf 0})-
G^{2}_{\lambda}({\bf 0},y)),$ $G_{\lambda}({\bf 0},r)/(\lambda+a)$
and their derivatives up to the order $[d/2]-1$ for each $r\in\mathbb{Z}^{d}$.

Integrating by parts it is easy to check the following equality
\begin{equation}\label{formula_for_K_d(lambda,y)}
K_{d}(\lambda;r)=\int\nolimits_{0}^{\infty}{e^{-\lambda
t}\left(\int\nolimits_{t}^{\infty}{p(u;{\bf
0},r)\,du}\right)\,dt},\quad r\in\mathbb{Z}^{d},\quad\lambda>0.
\end{equation}
Then for each $n\in\mathbb{N}$ function $K_{d}(\lambda;r)$ has the
$n$-th derivative w.r.t. $\lambda>0$ which can be calculated by the formula
\begin{equation}\label{derivative_K_d(lambda,y)}
K^{(n)}_{d}(\lambda;r)=(-1)^{n}\int\nolimits_{0}^{\infty}{e^{-\lambda
t}\,t^{n}\left(\int\nolimits_{t}^{\infty}{p(u;{\bf
0},r)\,du}\right)\,dt}
\end{equation}
for each $r\in\mathbb{Z}^{d}$. For convenience we set $K^{(0)}_{d}(\lambda;r):=K_{d}(\lambda;r)$.
Let us analyze for each $r$ the asymptotic behavior of the n-th derivative
($0\leq n\leq [d/2]-1$, $n\in\mathbb{Z}$) of the function $K_{d}(\lambda;r)$ at point $\lambda$, as $\lambda\to0+$.
To this end we need the relation implied by \eqref{p(t;x,y)sim} and Theorem 31 in \cite{V_Lectures},
namely,
\begin{equation}\label{int_p(t,x,y)sim}
\int\nolimits_{t}^{\infty}{p(u;{\bf
0},r)\,du}\sim\frac{2\gamma_{d}}{(d-2)t^{d/2-1}},\quad
t\to\infty,\quad r\in\mathbb{Z}^{d}.
\end{equation}
Then in view of \eqref{derivative_K_d(lambda,y)} one has, as $\lambda\to0+$,
\begin{equation}\label{derivative_K_d(lambda,y)_small}
(-1)^{n}K^{(n)}_{d}(\lambda;r)\to
\int\nolimits_{0}^{\infty}{t^{n}\left(\int\nolimits_{t}^{\infty}{p(u;{\bf
0},r)\,du}\right)\,dt}<\infty,\quad 0\leq
n<\left[\frac{d-1}{2}\right]-1.
\end{equation}
Moreover, using Tauberian theorem (Theorem 2 in \cite{Feller}, Ch.13, Sec.5) we observe that, as $\lambda\to0+$,
\begin{equation}\label{derivative_K_d(lambda,y)_small,only_even}
(-1)^{d/2-2}K^{(d/2-2)}_{d}(\lambda;r)\sim\frac{2\gamma_{d}}{(d-2)}\ln{\frac{1}{\lambda}}\quad\mbox{only for
even}\quad d.
\end{equation}
Note that for uneven $d$ the numbers $[(d-1)/2]-1$ and $[d/2]-1$ coincide whereas
for even $d$ one has $[(d-1)/2]-1=d/2-2$. Employing relations \eqref{derivative_K_d(lambda,y)} and
\eqref{int_p(t,x,y)sim} as well as Tauberian theorem (Theorem 4 in \cite{Feller}, Ch.13, Sec.5) we infer that,
as $\lambda\to0+$,
\begin{equation}\label{K([d/2]-1)_d(lambda;y)sim}
(-1)^{[d/2]-1}K_{d}^{([d/2]-1)}(\lambda;r)\sim\left\{
\begin{array}{lcl}
\frac{2\gamma_{d}\sqrt{\pi}}{(d-2)\sqrt{\lambda}},&\mbox{if}&\quad
d\quad\mbox{is uneven},\\
\frac{2\gamma_{d}}{(d-2)\lambda},&\mbox{if}&\quad d\quad\mbox{is even}.
\end{array}
\right.
\end{equation}
We stress that the right-hand sides of \eqref{derivative_K_d(lambda,y)_small,only_even} and
\eqref{K([d/2]-1)_d(lambda;y)sim} do not depend on point $r\in\mathbb{Z}^{d}$.

Now we investigate the asymptotic
properties of function $G_{\lambda}({\bf 0},{\bf0})-G_{\lambda}({\bf 0},y)$
and its derivatives up to the order $[d/2]-1$, as $\lambda\to0+$. Recall that for
$\lambda\geq0$ one has the identity
${G_{\lambda}({\bf 0},{\bf 0})-G_{\lambda}({\bf
0},y)=\int\nolimits_{0}^{\infty}{e^{-\lambda t}(p(t;{\bf 0},{\bf
0})-p(t;{\bf 0},y))\,dt}}$
where the asymptotic behavior of the function $p(t;{\bf 0},{\bf 0})-p(t;{\bf 0},y)$ (as $t\to\infty$)
is given by \eqref{p(t;0,0)-p(t,x,y)sim}. Then by virtue of the counterpart of \eqref{derivative_K_d(lambda,y)}
for the function $G_{\lambda}({\bf 0},{\bf 0})-G_{\lambda}({\bf 0},y)$ we get that, as $\lambda\to0+$,
\begin{equation}\label{derivative_G_lambda(0,0)-G_lambda(0,y)sim}
(-1)^{n}(G_{\lambda}({\bf 0},{\bf 0})-G_{\lambda}({\bf
0},y))^{(n)}\to\int\nolimits_{0}^{\infty}{t^{n}(p(t;{\bf 0},{\bf
0})-p(t;{\bf 0},y))\,dt}<\infty,\quad 0\leq
n\leq\left[\frac{d}{2}\right]-1.
\end{equation}

Using Theorem \ref{Faa_di_Bruno} we are ready to estimate the asymptotic growth of the $n$-th
derivative, ${0\leq n\leq[d/2]-1}$, of the function $(G_{\lambda}({\bf 0},{\bf 0})-G_{\lambda}({\bf 0},y))^{-1}$,
as $\lambda\to0+$. For this purpose we consider function $V^{-1}$ as the external function $W(V)$
and we substitute $G_{\lambda}({\bf 0},{\bf 0})-G_{\lambda}({\bf 0},y)$ instead of internal function
$V(\cdot)$ in Theorem \ref{Faa_di_Bruno}. Since $W^{(n)}(V)=\left(V^{-1}\right)^{(n)}=(-1)^{n}n!\,V^{-n-1}$
for every $n\in\mathbb{Z}_{+}$ and $G_{\lambda}({\bf 0},{\bf 0})-G_{\lambda}({\bf 0},y)\to a^{-1}\rho_{d}(y)>0$,
as $\lambda\to0+$, then in view of Theorem \ref{Faa_di_Bruno} and relation
\eqref{derivative_G_lambda(0,0)-G_lambda(0,y)sim} we obtain the estimate we are interested in
\begin{equation}\label{derivative_(G_lambda(0,0)-G_lambda(0,y))^-1}
\left(\frac{1}{G_{\lambda}({\bf 0},{\bf 0})-G_{\lambda}({\bf
0},y)}\right)^{(n)}=O(1),\quad\lambda\to0+,\quad 0\leq
n\leq\left[\frac{d}{2}\right]-1.
\end{equation}

Let us turn to the study of the asymptotic behavior of the function $G_{\lambda}({\bf 0},r)$ and its $n$-th
derivatives, $0<n\leq [d/2]-1$, as $\lambda\to0+$, for each
$r\in\mathbb{Z}^{d}$. In view of the counterpart of \eqref{derivative_K_d(lambda,y)} for
$G_{\lambda}({\bf 0},r)$ as well as of relation \eqref{p(t;x,y)sim} we see that, as $\lambda\to0+$,
\begin{equation}\label{derivative_G_lambda(0,y)sim}
(-1)^{n}G^{(n)}_{\lambda}({\bf
0},r)\to\int\nolimits_{0}^{\infty}{t^{n}p(t;{\bf
0},r)dt}<\infty,\quad 0\leq n\leq\left[\frac{d-1}{2}\right]-1.
\end{equation}
Moreover, employing Tauberian theorem (Theorem 2 in \cite{Feller}, Ch.13, Sec.5)
we find that
\begin{equation}\label{derivative_G_lambda(0,y)sim,only_even}
(-1)^{d/2-1}G_{\lambda}^{(d/2-1)}({\bf
0},r)\sim\gamma_{d}\ln{\frac{1}{\lambda}},\quad\lambda\to0+,\quad\mbox{only for even}\quad d.
\end{equation}
As mentioned above, the numbers $[(d-1)/2]-1$ and $[d/2]-1$ coincide for uneven $d$
whereas for even $d$ one has $[(d-1)/2]-1=d/2-2$.

Now for estimating the asymptotic behavior of the $n$-th order derivatives,
$0\leq n\leq[d/2]-1$, of the function $G_{\lambda}({\bf 0},r)/(\lambda+a)$,
as $\lambda\to0+$, it suffices to apply the Leibniz formula
as well as relations \eqref{derivative_G_lambda(0,y)sim} and \eqref{derivative_G_lambda(0,y)sim,only_even}.
Indeed,
$$\left(\frac{G_{\lambda}({\bf 0},r)}{\lambda+a}\right)^{(n)}=\sum\limits_{k=0}^{n}{C^{k}_{n}\,
G^{(k)}_{\lambda}({\bf 0},r)\left(\frac{1}{\lambda+a}\right)^{(n-k)}}=\sum\limits_{k=0}^{n}{C^{k}_{n}\,
G^{(k)}_{\lambda}({\bf 0},r)\frac{(-1)^{n-k}(n-k)!\,}{(\lambda+a)^{n-k+1}}}$$
where $n\in\mathbb{Z}_{+}$ and, hence, for each $r\in\mathbb{Z}^{d}$ we obtain the desired estimates
\begin{equation}\label{derivative_G_lambda(0,y)/(lambda+a)sim}
\left(\frac{G_{\lambda}({\bf
0},r)}{\lambda+a}\right)^{(n)}=O(1),\quad\lambda\to0+,\quad0\leq
n\leq\left[\frac{d-1}{2}\right]-1,
\end{equation}
\begin{equation}\label{derivative_G_lambda(0,y)/(lambda+a)sim_only_even}
\left(\frac{G_{\lambda}({\bf
0},r)}{\lambda+a}\right)^{(d/2-1)}=O\left(\ln{\frac{1}{\lambda}}\right),\quad\lambda\to0+,\quad\mbox{only for
even}\quad d.
\end{equation}

Repeating the argument used for deriving relations \eqref{derivative_G_lambda(0,y)sim} and
\eqref{derivative_G_lambda(0,y)sim,only_even} we find the asymptotic behavior of the function
$G_{\lambda}({\bf 0},{\bf 0})+G_{\lambda}({\bf 0},y)$ and its derivatives up to the order $[d/2]-1$,
as $\lambda\to0+$,
\begin{equation}\label{derivative_G_lambda(0,0)+G_lambda(0,y)sim}
(-1)^{n}(G_{\lambda}({\bf 0},{\bf 0})+G_{\lambda}({\bf
0},y))^{(n)}\to\int\limits_{0}^{\infty}{t^{n}(p(t;{\bf 0},{\bf
0})+p(t;{\bf 0},y))dt}<\infty,\,0\leq
n\leq\left[\frac{d-1}{2}\right]-1,
\end{equation}
\begin{equation}\label{derivative_G_lambda(0,0)+G_lambda(0,y)sim,only_even}
(-1)^{d/2-1}\left(G_{\lambda}({\bf 0},{\bf 0})+G_{\lambda}({\bf
0},y)\right)^{(d/2-1)}\sim2\gamma_{d}\ln{\frac{1}{\lambda}}\quad\mbox{only for even}\quad d.
\end{equation}

To investigate the asymptotic properties of the function $(G_{\lambda}({\bf 0},{\bf 0})+G_{\lambda}({\bf 0},y))^{-1}$
and its derivatives of order $n$, $0<n\leq[d/2]-1$, as $\lambda\to0+$, we employ Theorem \ref{Faa_di_Bruno}
once again. We take $V^{-1}$ as the external function $W(V)$ and consider the function
$G_{\lambda}({\bf 0},{\bf 0})+G_{\lambda}({\bf 0},y)$ as $V(\cdot)$ in Theorem~\ref{Faa_di_Bruno}. Since
$W^{(n)}(V)=\left(V^{-1}\right)^{(n)}=(-1)^{n}n!\,V^{-n-1}$ for all $n\in\mathbb{Z}_{+}$
and ${G_{\lambda}({\bf 0},{\bf 0})+G_{\lambda}({\bf 0},y)\to G_{0}({\bf 0},{\bf 0})+G_{0}({\bf
0},y)>0}$, as $\lambda\to0+$, due to Theorem \ref{Faa_di_Bruno} as well as relations
\eqref{derivative_G_lambda(0,0)+G_lambda(0,y)sim} and
\eqref{derivative_G_lambda(0,0)+G_lambda(0,y)sim,only_even} we get the estimates
\begin{equation}\label{derivative_(G_lambda(0,0)+G_lambda(0,0))^-1}
\left(\frac{1}{G_{\lambda}({\bf 0},{\bf 0})+G_{\lambda}({\bf
0},y)}\right)^{(n)}=O(1),\quad \lambda\to0+,\quad 0\leq
n\leq\left[\frac{d-1}{2}\right]-1,
\end{equation}
\begin{equation}\label{derivative_(G_lambda(0,0)+G_lambda(0,0))^-1_only_even}
\left(\frac{1}{G_{\lambda}({\bf 0},{\bf 0})+G_{\lambda}({\bf
0},y)}\right)^{(d/2-1)}=O\left(\ln{\frac{1}{\lambda}}\right),\quad\lambda\to0+,\quad\mbox{only for even}\quad d.
\end{equation}

To establish asymptotic behavior of the function
$(G^{2}_{\lambda}({\bf 0},{\bf 0})-G^{2}_{\lambda}({\bf 0},y))^{-1}$ and its derivatives up to the order
$[d/2]-1$, as $\lambda\to0+$, it suffices to apply the Leibniz formula along with estimates
\eqref{derivative_(G_lambda(0,0)-G_lambda(0,y))^-1},
\eqref{derivative_(G_lambda(0,0)+G_lambda(0,0))^-1} and
\eqref{derivative_(G_lambda(0,0)+G_lambda(0,0))^-1_only_even}.
Indeed,
$$\left(\frac{1}{G^{2}_{\lambda}({\bf 0},{\bf 0})-G^{2}_{\lambda}({\bf 0},y)}\right)^{(n)}
=\sum\limits_{k=0}^{n}{C^{k}_{n}\left(\frac{1}{G_{\lambda}({\bf 0},{\bf 0})
-G_{\lambda}({\bf 0},y)}\right)^{(k)}\left(\frac{1}{G_{\lambda}({\bf 0},{\bf 0})
+G_{\lambda}({\bf 0},y)}\right)^{(n-k)}}$$
for $n\in\mathbb{Z}_{+}$, and, consequently,
\begin{equation}\label{derivative_(G^2_lambda(0,0)-G^2_lambda(0,y))^-1}
\left(\frac{1}{G^{2}_{\lambda}({\bf 0},{\bf 0})-G^{2}_{\lambda}({\bf
0},y)}\right)^{(n)}=O(1),\quad\lambda\to0+,\quad 0\leq
n\leq\left[\frac{d-1}{2}\right]-1,
\end{equation}
\begin{equation}\label{derivative_(G^2_lambda(0,0)-G^2_lambda(0,y))^-1_only_even}
\left(\frac{1}{G^{2}_{\lambda}({\bf 0},{\bf 0})-G^{2}_{\lambda}({\bf
0},y)}\right)^{(d/2-1)}=O\left(\ln{\frac{1}{\lambda}}\right),\quad\lambda\to0+,\quad\mbox{only for even}\quad d.
\end{equation}

Let us investigate the asymptotic properties of the function
${G_{\lambda}({\bf 0},r)/(G^{2}_{\lambda}({\bf 0},{\bf
0})-G^{2}_{\lambda}({\bf 0},y))}$ and its $n$-th derivatives,
$0\leq n\leq[d/2]-1$, for each $r\in\mathbb{Z}^{d}$, as $\lambda\to0+$. The Leibniz formula implies
that for all $\lambda>0$ and $n\in\mathbb{Z}_{+}$
$$\left(\frac{G_{\lambda}({\bf 0},r)}{G^{2}_{\lambda}({\bf 0},{\bf 0})-G^{2}_{\lambda}({\bf 0},y)}\right)^{(n)}
=\sum\limits_{k=0}^{n}{C^{k}_{n}G^{(k)}_{\lambda}({\bf 0},r)\left(\frac{1}{G^{2}_{\lambda}({\bf 0},{\bf 0})
-G^{2}_{\lambda}({\bf 0},y)}\right)^{(n-k)}}.$$
Taking into account \eqref{derivative_G_lambda(0,y)sim},
\eqref{derivative_G_lambda(0,y)sim,only_even},
\eqref{derivative_(G^2_lambda(0,0)-G^2_lambda(0,y))^-1} and
\eqref{derivative_(G^2_lambda(0,0)-G^2_lambda(0,y))^-1_only_even},
for each $r\in\mathbb{Z}^{d}$, we infer that
\begin{equation}\label{derivative_G_lambda(0,z)(G^2_lambda(0,0)-G^2_lambda(0,y))^-1}
\left(\frac{G_{\lambda}({\bf 0},r)}{G^{2}_{\lambda}({\bf 0},{\bf
0})-G^{2}_{\lambda}({\bf
0},y)}\right)^{(n)}=O(1),\quad\lambda\to0+,\quad0\leq
n\leq\left[\frac{d-1}{2}\right]-1,
\end{equation}
\begin{equation}\label{derivative_G_lambda(0,z)(G^2_lambda(0,0)-G^2_lambda(0,y))^-1_only_even}
\left(\frac{G_{\lambda}({\bf 0},r)}{G^{2}_{\lambda}({\bf 0},{\bf
0})-G^{2}_{\lambda}({\bf
0},y)}\right)^{(d/2-1)}=O\left(\ln{\frac{1}{\lambda}}\right),\quad\lambda\to0+,\quad\mbox{only for even}\quad d.
\end{equation}

Since all the basic relations are established, we turn to finding the asymptotic behavior of the function
$J^{([d/2]-1)}_{1}(\lambda;r,y)$ for each $r\in\mathbb{Z}^{d}$, as $\lambda\to0+$.
By the definition of $J_{1}(\lambda;r,y)$ after applying the Leibniz formula to its
$([d/2]-1)$-th derivative we obtain
$$J^{([d/2]-1)}_{1}(\lambda;r,y)=\sum\limits_{k=0}^{[d/2]-1}{C^{k}_{[d/2]-1}\left(\frac{1}
{G^{2}_{\lambda}({\bf 0},{\bf 0})-G^{2}_{\lambda}({\bf 0},y)}\right)^{(k)}K^{([d/2]-1-k)}_{d}(\lambda;r)}.$$
Analyzing every summand in the last sum with the help of relations
\eqref{derivative_K_d(lambda,y)_small}--\eqref{K([d/2]-1)_d(lambda;y)sim} and
\eqref{derivative_(G^2_lambda(0,0)-G^2_lambda(0,y))^-1},
\eqref{derivative_(G^2_lambda(0,0)-G^2_lambda(0,y))^-1_only_even}
we conclude that the first summand gives the main contribution to the asymptotic behavior of
$J^{([d/2]-1)}_{1}(\lambda;x,y)$, $\lambda\to0+$. Indeed, the first summand has the order of growth
$1/\sqrt{\lambda}$ for uneven $d$ and has the order $1/\lambda$ for even $d$, respectively.
On the other hand, the rest summands are $O(1)$ for uneven $d$ and are $O(\ln^{2}{\lambda})$ for even $d$,
respectively. Therefore, for each $r\in\mathbb{Z}^{d}$, as $\lambda\to0+$, one gets
\begin{equation}\label{J^([d/2]-1)_1sim}
(-1)^{\left[\frac{d}{2}\right]-1}J^{\left(\left[\frac{d}{2}\right]-1\right)}_{1}(\lambda;r,y)\sim\left\{
\begin{array}{lcl}
\frac{2\gamma_{d}\sqrt{\pi}}{(d-2)(G^{2}_{0}({\bf 0},{\bf
0})-G^{2}_{0}({\bf 0},y))\sqrt{\lambda}},&\mbox{if}&
d\quad\mbox{is uneven},\\
\frac{2\gamma_{d}}{(d-2)(G^{2}_{0}({\bf 0},{\bf 0})-G^{2}_{0}({\bf
0},y))\lambda},&\mbox{if}& d\quad\mbox{is even}.
\end{array}
\right.
\end{equation}

To estimate the asymptotic behavior of the functions $J^{([d/2]-1)}_{2}(\lambda;r,y)$ and
$J^{([d/2]-1)}_{3}(\lambda;r,y)$, as $\lambda\to0+$, we use the Leibniz formula once again, namely,
$$J^{([d/2]-1)}_{2}(\lambda;r,y)=\sum\limits_{k=0}^{[d/2]-1}{C^{k}_{[d/2]-1}\left(\frac{G_{\lambda}({\bf
0},{\bf 0})}{G^{2}_{\lambda}({\bf 0},{\bf 0})-G^{2}_{\lambda}({\bf
0},y)}\right)^{(k)}\left(K_{d}(\lambda;y-r)-K_{d}(\lambda;y)\right)^{([d/2]-1-k)}},$$
$$J^{([d/2]-1)}_{3}(\lambda;r,y)=\sum\limits_{k=0}^{[d/2]-1}{C^{k}_{[d/2]-1}\left(\frac{G_{\lambda}({\bf
0},y)}{G^{2}_{\lambda}({\bf 0},{\bf 0})-G^{2}_{\lambda}({\bf
0},y)}\right)^{(k)}\left(K_{d}(\lambda;{\bf
0})-K_{d}(\lambda;r)\right)^{([d/2]-1-k)}}.$$
By virtue of \eqref{derivative_K_d(lambda,y)_small}--\eqref{K([d/2]-1)_d(lambda;y)sim},
\eqref{derivative_G_lambda(0,z)(G^2_lambda(0,0)-G^2_lambda(0,y))^-1}
and \eqref{derivative_G_lambda(0,z)(G^2_lambda(0,0)-G^2_lambda(0,y))^-1_only_even} at the right-hand sides
of the last equalities there are only the first summands which grow as $o(1/\sqrt{\lambda})$ for uneven $d$
and as $o(1/\lambda)$ for even $d$, respectively, as $\lambda\to0+$. The rest summands are $O(1)$
for uneven $d$ and are $o(\ln^{2}{\lambda})$ for even $d$. Thus, on account of \eqref{J^([d/2]-1)_1sim}
for each $r\in\mathbb{Z}^{d}$, as $\lambda\to0+$,
we come to the desired estimates
\begin{equation}\label{J^([d/2]-1)_5sim,J^([d/2]-1)_5sim}
J^{([d/2]-1)}_{2}(\lambda;r,y)=o\left(J^{([d/2]-1)}_{1}(\lambda;{\bf
0},y)\right),\quad
J^{([d/2]-1)}_{3}(\lambda;r,y)=o\left(J^{([d/2]-1)}_{1}(\lambda;{\bf
0},y)\right).
\end{equation}

Let us study the asymptotic growth of the function $J^{([d/2]-1)}_{4}(\lambda;y)$ as $\lambda\to0+$.
The definition of $J_{4}(\lambda;y)$ and the Leibniz formula entail
$$J^{([d/2]-1)}_{4}(\lambda;y)=\sum\limits_{k=0}^{[d/2]-1}{C^{k}_{[d/2]-1}\left(\frac{1}
{G^{2}_{\lambda}({\bf 0},{\bf 0})-G^{2}_{\lambda}({\bf 0},y)}\right)^{(k)}
\left(K_{d}(\lambda;{\bf 0})-K_{d}(\lambda;y)\right)^{([d/2]-1-k)}}.$$
Due to relations \eqref{derivative_K_d(lambda,y)_small}--\eqref{K([d/2]-1)_d(lambda;y)sim},
\eqref{derivative_(G^2_lambda(0,0)-G^2_lambda(0,y))^-1} and
\eqref{derivative_(G^2_lambda(0,0)-G^2_lambda(0,y))^-1_only_even} we note that
the first summand at the last sum is $o(1/\sqrt{\lambda})$ for uneven $d$ and is $o(1/\lambda)$ for even $d$,
respectively, as $\lambda\to0+$. Meanwhile, the rest summands have the order of growth $O(1)$ for uneven $d$
and $o(\ln^{2}{\lambda})$ for even $d$, respectively. Hence, taking into account
\eqref{J^([d/2]-1)_1sim} we conclude that for each $r\in\mathbb{Z}^{d}$
\begin{equation}\label{J^([d/2]-1)_2sim}
J^{([d/2]-1)}_{4}(\lambda;y)=o\left(J^{([d/2]-1)}_{1}(\lambda;r,y)\right),\quad\lambda\to0+.
\end{equation}

We turn to estimating of the asymptotic behavior of the function $J^{([d/2]-1)}_{5}(\lambda;y)$,
as $\lambda\to0+$. The definition of $J_{5}(\lambda;y)$ and the Leibniz formula imply that
$$J^{([d/2]-1)}_{5}(\lambda;y)=\sum\limits_{k=0}^{[d/2]-1}{C^{k}_{[d/2]-1}\left(\frac{1}{G_{\lambda}
({\bf 0},{\bf 0})-G_{\lambda}({\bf 0},y)}\right)^{(k)}\left(K_{d}(\lambda;{\bf 0})-
K_{d}(\lambda;y)\right)^{([d/2]-1-k)}}.$$
On account of relations \eqref{derivative_K_d(lambda,y)_small}--\eqref{K([d/2]-1)_d(lambda;y)sim}
and \eqref{derivative_(G_lambda(0,0)-G_lambda(0,y))^-1} we see that
the first summand at the last sum is $o(1/\sqrt{\lambda})$ for uneven $d$ and is $o(1/\lambda)$
for even $d$, respectively. All the other summands have the order of growth $O(1)$ for uneven $d$ and
$o(\ln{\lambda})$ for even $d$, respectively, as $\lambda\to0+$. By \eqref{J^([d/2]-1)_1sim}
it follows that for each $r\in\mathbb{Z}^{d}$
\begin{equation}\label{J^([d/2]-1)_3sim}
J^{([d/2]-1)}_{5}(\lambda;y)=o\left(J^{([d/2]-1)}_{1}(\lambda;r,y)\right),\quad\lambda\to0+.
\end{equation}

The last step is the revealing the asymptotic properties, as $\lambda\to0+$, of the function
$J^{([d/2]-1)}_{6}(\lambda;r,y)$. By the definition of $J_{6}(\lambda;r,y)$ together with the Leibniz formula
one can write
$$J^{([d/2]-1)}_{6}(\lambda;r,y)=\sum\limits_{k=0}^{[d/2]-1}{C^{k}_{[d/2]-1}\left(\frac{1}{G^{2}_{\lambda}
({\bf 0},{\bf 0})-G^{2}_{\lambda}({\bf 0},y)}\right)^{(k)}\left(\frac{G_{\lambda}({\bf 0},r)}
{\lambda+a}\right)^{([d/2]-1-k)}}.$$
Recalling relations \eqref{derivative_G_lambda(0,y)/(lambda+a)sim},
\eqref{derivative_G_lambda(0,y)/(lambda+a)sim_only_even},
\eqref{derivative_(G^2_lambda(0,0)-G^2_lambda(0,y))^-1} and
\eqref{derivative_(G^2_lambda(0,0)-G^2_lambda(0,y))^-1_only_even}, we deduce that if
$d$ is uneven then all the summands at the last sum are $O(1)$, as $\lambda\to0+$.
However, if $d$ is even then the first and the last summands have the order of growth $O(\ln{\lambda})$
whereas the others have the order of growth $O(1)$, as $\lambda\to0+$.
In any case in view of \eqref{J^([d/2]-1)_1sim} the following formula holds true
\begin{equation}\label{J^([d/2]-1)_4sim}
J^{([d/2]-1)}_{6}(\lambda;r,y)=o\left(J^{([d/2]-1)}_{1}(\lambda;r,y)\right),\quad\lambda\to0+.
\end{equation}

Combining results \eqref{lt_Hxy0(infty)-Hxy0_d>=3_sum} -- \eqref{lt_Hooy(infty)-Hooy_d>=3_sum} and
\eqref{J^([d/2]-1)_1sim} -- \eqref{J^([d/2]-1)_4sim} we conclude that
\begin{eqnarray}\label{Hxy0(infty)-Hxy0(lambda)_derivative}
\int\nolimits_{0}^{\infty}{e^{-\lambda t}t^{[d/2]-1}(H_{x,y,{\bf
0}}(\infty)-H_{x,y,{\bf
0}}(t))\,dt}&\sim&\frac{2\gamma_{d}(\rho_{d}(y-x)+\rho_{d}(x)-\rho_{d}(y))}{a(d-2)(G_{0}({\bf
0},{\bf 0})+G_{0}({\bf
0},y))^{2}}\nonumber\\
&\times&\left\{
\begin{array}{lcl}
\sqrt{\pi}/\sqrt{\lambda},&\mbox{if}&
d\quad\mbox{is uneven},\\
1/\lambda,&\mbox{if}& d\quad\mbox{is even},
\end{array}
\right.
\end{eqnarray}
\begin{eqnarray}\label{H0y0(infty)-H0y0(lambda)_derivative}
\int\nolimits_{0}^{\infty}{e^{-\lambda
t}t^{\left[\frac{d}{2}\right]-1}(H_{{\bf 0},y,{\bf
0}}(\infty)-H_{{\bf 0},y,{\bf 0}}(t))\,dt}\sim\int\nolimits_{0}^{\infty}{e^{-\lambda
t}t^{\left[\frac{d}{2}\right]-1}(H_{y,y,{\bf 0}}(\infty)-H_{y,y,{\bf 0}}(t))\,dt}\nonumber\\
\sim\frac{2\gamma_{d}}{a(d-2)(G_{0}({\bf 0},{\bf 0})+G_{0}({\bf
0},y))^{2}}\left\{
\begin{array}{lcl}
\sqrt{\pi}/\sqrt{\lambda},&\mbox{if}&
d\quad\mbox{is uneven},\\
1/\lambda,&\mbox{if}& d\quad\mbox{is even},
\end{array}
\right.
\end{eqnarray}
as $\lambda\to0+$. If $C_{d}(x,y)>0$ for all $x,y\in\mathbb{Z}^{d}\backslash\{{\bf 0}\},$ $x\neq y$
(it is easily seen that $C_{d}({\bf 0},y)$ and $C_{d}(y,y)$ are always positive for $d\geq3$) then
application of Corollary 43 in \cite{V_Lectures} to the last asymptotic relations
implies \eqref{H_xyz(infty)-H_xyz(t)sim,d>=3}.

Thus, to complete the proof of Theorem \ref{T:not_simple} for $d\geq3$ we only have to verify the positivity
of $C_{d}(x,y)$ for $x,y\in\mathbb{Z}^{d}\backslash\{{\bf 0}\}$, $x\neq y$.
The main idea of the demonstration is the same as that exploited for proving
the positivity of $C_{2}(\cdot,\cdot)$ in Subsection 4.2. Namely, we assume the contrary, that is,
for some $d\in\mathbb{N}$, $x\in\mathbb{Z}^{d}\backslash\{{\bf 0},y\}$
and $y\in\mathbb{Z}^{d}\backslash\{{\bf 0}\}$ the relation $C_{d}(x,y)=0$ is valid.
This is equivalent to the following relations
\begin{equation}\label{Cd(x,y)=0...,d>=3}
G_{0}({\bf 0},{\bf 0})+G_{0}({\bf 0},y)=G_{0}(x,y)+G_{0}({\bf
0},x)\Leftrightarrow G_{0}({\bf 0},{\bf 0})-G_{0}({\bf
0},x)=G_{0}(x,y)-G_{0}({\bf 0},y).
\end{equation}
Note that the formula for $H_{y,x,{\bf 0}}(\infty)$ appearing in \eqref{T_Hxyz(infty),d>=3} and the definition
of function $\rho_{d}(\cdot)$ allow us to write
$$H_{y,x,{\bf 0}}(\infty)=\frac{G_{0}({\bf 0},{\bf
0})(G_{0}(x,y)-G_{0}({\bf 0},y))+G_{0}({\bf 0},y)(G_{0}({\bf 0},{\bf
0})-G_{0}({\bf 0},x))}{G^{2}_{0}({\bf 0},{\bf 0})-G^{2}_{0}({\bf
0},x)}.$$
Combining \eqref{Cd(x,y)=0...,d>=3} and the last equality we see that
\begin{eqnarray*}
H_{y,x,{\bf 0}}(\infty)&=&\frac{(G_{0}({\bf 0},{\bf 0})-G_{0}({\bf 0},x))(G_{0}({\bf
0},{\bf 0})+G_{0}({\bf 0},y))}{G^{2}_{0}({\bf 0},{\bf
0})-G^{2}_{0}({\bf 0},x)}\\
&=&\frac{G_{0}({\bf 0},{\bf 0})+G_{0}({\bf
0},y)}{G_{0}({\bf 0},{\bf 0})+G_{0}({\bf
0},x)}=\frac{G_{0}(x,y)+G_{0}({\bf 0},x)}{G_{0}({\bf 0},{\bf
0})+G_{0}({\bf 0},x)}.
\end{eqnarray*}
However, by formula \eqref{main_equation} one has
$H_{y,x,{\bf 0}}(\infty)\leq H_{y,x}(\infty)$ and due to \eqref{Hxy(infty),d>=3} it follows that
$$\frac{G_{0}(x,y)+G_{0}({\bf 0},x)}{G_{0}({\bf 0},{\bf
0})+G_{0}({\bf 0},x)}\leq\frac{G_{0}(x,y)}{G_{0}({\bf 0},{\bf
0})}\Leftrightarrow G_{0}({\bf 0},{\bf 0})\leq G_{0}(x,y).$$
The obtained contradiction completes the proof of \eqref{H_xyz(infty)-H_xyz(t)sim,d>=3} and the proof of
Theorem \ref{T:not_simple} as well.

\section{Proof of Theorem \ref{T:simple}}

Firstly note that in view of Corollary \ref{C:properties_of_Hxyz} it suffices to establish Theorem
\ref{T:simple} for $z=0$. Recall that for a simple random walk on $\mathbb{Z}$ the
transition intensities $a(r,r+1)=a(r,r-1)=a/2$ and $a(r,r+k)=0$ for all $r\in\mathbb{Z}$ and
$k\in\mathbb{Z}$ such that $|k|>1$. By definition of the function $\phi(\theta)$ and by that
of the constant $\gamma_{1}$
one gets $\phi(\theta)=a(\cos{\theta}-1)$, $\theta\in[-\pi,\pi]$, and $\gamma_{1}=1/\sqrt{2 a \pi}$.
Hence, equality \eqref{rho_d(x)=int}
and Lemma \ref{L:trigonometric_equality} imply that $\rho_{1}(r)=|r|$, $r\in\mathbb{Z}\backslash\{{\bf 0}\}$.
Observe that the formula for $H_{x,y,{\bf 0}}(\infty)$ in \eqref{T_Hxyz(infty),d=1,2} is valid
for a simple random walk on $\mathbb{Z}$ as well. Therefore, we infer that
$H_{x,y,{\bf 0}}(\infty)=0$ for $y<{\bf 0}<x$ and $x<{\bf 0}<y$, $H_{x,y,{\bf 0}}(\infty)=x/y$ for
${\bf 0}<x<y$ and $y<x<{\bf 0}$, $H_{x,y,{\bf 0}}(\infty)=1$ for ${\bf 0}<y<x$ and $x<y<{\bf 0}$,
$H_{{\bf 0},y,{\bf 0}}(\infty)=1/(2|y|)$ and $H_{y,y,{\bf 0}}(\infty)=1-1/(2|y|)$.
Thus, for $y<{\bf 0}<x$ and $x<{\bf 0}<y$ one has $H_{x,y,{\bf 0}}(\cdot)\equiv0$, that is, relation
\eqref{Hxyz_simple_5} is proved. For ${\bf 0}<y<x$ and $x<y<{\bf 0}$ due to \eqref{main_equation}
we conclude that $H_{x,y,{\bf0}}(\infty)=1=H_{x,y}(\infty)$ and $H_{x,{\bf 0},y}(\infty)=0$.
This is equivalent to relations $H_{x,{\bf 0},y}(\cdot)\equiv0$ and $H_{x,y,{\bf 0}}(\cdot)\equiv H_{x,y}(\cdot)$.
Since the asymptotic behavior of the function $H_{x,y}(\infty)-H_{x,y}(t)$ is established in
Lemma \ref{L:Hxy}, relation \eqref{Hxyz_simple_1} is also proved. Moreover,
formula \eqref{lt_Hxy0(infty)-Hxy0_sim} is true for a simple random walk on $\mathbb{Z}$ and,
consequently, applying Tauberian theorem (Theorem 4 in \cite{Feller}, Ch.13, Sec.5) we find the asymptotic
behavior of $H_{y,y,{\bf 0}}(\infty)-H_{y,y,{\bf 0}}(t)$, as $t\to\infty$, with $C_{1}(y,y)=1/(2\sqrt{a \pi})>0$.
Thus, relation \eqref{Hxyz_simple_2} is proved. So we have to investigate only the asymptotic properties
of $H_{x,y,{\bf 0}}(\infty)-H_{x,y,{\bf 0}}(t)$ when $x\in[0,y),$ $x\in\mathbb{Z}$, that is,
to verify relations \eqref{Hxyz_simple_3} and \eqref{Hxyz_simple_4}.

To this end recall a well-known result (the gambler ruin problem) for embedded chain
$\{S_{n}, n\in\mathbb{Z}_{+}\}$. Namely, ${\sf P}({\bf 0}<S_{k}<y, 0<k\leq n|S(0)=x)\leq (1-
\varepsilon_{0})^{n}$ for some $\varepsilon_{0}\in(0,1)$ and $x\in[0,y],$
$x\in\mathbb{Z}$ (see, e.g., \cite{Shiriaev}, Ch.1, Sec.9).
Bearing on this result we derive that
\begin{eqnarray*}
H_{x,y,{\bf 0}}(\infty)-H_{x,y,{\bf
0}}(t)&=&\sum\nolimits_{n=0}^{\infty}{{\sf P}\left(\left.t<\tau_{y,{\bf
0}}<\infty, N(t)=n\right|S(0)=x\right)}\\
&\leq&\sum\nolimits_{n=0}^{\infty}{{\sf P}\left(\left.{\bf
0}<S_{k}<y,0<k\leq n, N(t)=n\right|S(0)=x\right)}\\
&=&\sum\nolimits_{n=0}^{\infty}{{\sf P}\left(\left.{\bf
0}<S_{k}<y,0<k\leq n\right|S(0)=x\right)\,{\sf
P}\left(\left.N(t)=n\right|S(0)=x\right)}\\
&\leq&\sum\nolimits_{n=0}^{\infty}{(1-\varepsilon_{0})^{n}\frac{(a t)^{n}e^{-a
t}}{n!}}=e^{-a\,\varepsilon_{0}\,t}.
\end{eqnarray*}
It follows that relations \eqref{Hxyz_simple_3} and \eqref{Hxyz_simple_4} are valid with
$\varepsilon\in(0,\varepsilon_{0})$. Theorem \ref{T:simple} is proved completely.

\section{Proof of Theorem \ref{T:H-}}

On account of Corollary \ref{C:properties_of_Hxyz} it is sufficient to prove Theorem \ref{T:H-} for $z=0$.
Since ${\tau_{y,{\bf 0}}=\tau^{-}_{y,{\bf 0}}+\tau}$ a.s., one has $H_{x,y,{\bf 0}}(t)=H^{-}_{x,y,{\bf 0}}\ast G(t)$,
$t\geq0$. This immediately implies ${H_{x,y,{\bf 0}}(\infty)=H^{-}_{x,y,{\bf 0}}(\infty)}$,
$x\in\mathbb{Z}^{d}$, $y\in\mathbb{Z}^{d}\backslash\{{\bf 0}\}$.
Furthermore, by formula \eqref{universal_formula} we get that for $\lambda>0$
\begin{equation}\label{lt_H-(infty)-H-}
\widehat{(H^{-}_{x,y,{\bf 0}}(\infty)-H^{-}_{x,y,{\bf 0}})}(\lambda)=\widehat{(H_{x,y,{\bf 0}}(\infty)-
H_{x,y,{\bf 0}})}(\lambda)-a^{-1}\check{H}_{x,y,{\bf 0}}(\lambda).
\end{equation}
In view of \eqref{H_xyz(infty)-H_xyz(t)sim,d=1} -- \eqref{Hxyz_simple_2}
and Tauberian theorem (Theorem 4 in \cite{Feller}, Ch.13, Sec.5)
it is easy to see that the first summand at the right-hand side of \eqref{lt_H-(infty)-H-}
makes the main contribution to the asymptotic behavior of the left-hand side of \eqref{lt_H-(infty)-H-}.
More exactly, this is true if $d\leq3$, $x\in\mathbb{Z}^{d}$, $y\in\mathbb{Z}^{d}\backslash\{{\bf 0}\}$,
except for a simple random walk on $\mathbb{Z}$ in the cases
${\bf 0}\leq x<y$, $y<x\leq{\bf 0}$, $x<{\bf 0}<y$ or $y<{\bf 0}<x$.
Therefore, $\widehat{(H^{-}_{x,y,{\bf 0}}(\infty)-
H^{-}_{x,y,{\bf 0}})}(\lambda)\sim\widehat{(H_{x,y,{\bf 0}}(\infty)-H_{x,y,{\bf 0}})}(\lambda)$, as
$\lambda\to0+$, and due to \eqref{H_xyz(infty)-H_xyz(t)sim,d=1}--\eqref{Hxyz_simple_2}
and Tauberian theorem (Theorem 4 in \cite{Feller}, Ch.13, Sec.5)
we prove Theorem \ref{T:H-} for the mentioned $d$, $x$ and $y$.

When the random walk on $\mathbb{Z}$ is simple and ${\bf 0}\leq x<y$ or $y<x\leq{\bf 0}$,
the estimate of $H^{-}_{x,y,{\bf 0}}(\infty)-H^{-}_{x,y,{\bf 0}}(t)$, as $t\to\infty$, is a direct consequence of Theorem \ref{T:simple} and
the inequality $H_{x,y,{\bf 0}}(\infty)-H_{x,y,{\bf 0}}(t)\geq H^{-}_{x,y,{\bf 0}}(\infty)-H^{-}_{x,y,{\bf 0}}(t)$
used in the proof of Theorem \ref{T:not_simple} for $d=1$. Moreover, for a simple random walk on $\mathbb{Z}$
when $x<{\bf 0}<y$ or $y<{\bf 0}<x$ we may assert that $H^{-}_{x,y,{\bf 0}}(t)\equiv 0$ by virtue of
the established identity $H^{-}_{x,y,{\bf 0}}(\infty)=H_{x,y,{\bf 0}}(\infty)$ and relation \eqref{Hxyz_simple_5}.
Thus, for these cases Theorem \ref{T:H-} is also proved.

To find the asymptotic behavior of $H^{-}_{x,y,{\bf 0}}(\infty)-H^{-}_{x,y,{\bf 0}}(t)$, as $t\to\infty$,
for $d\geq4$, we differentiate $[d/2]-1$ times the both sides of equality \eqref{lt_H-(infty)-H-}
w.r.t. $\lambda>0$.
The asymptotic relation for the $([d/2]-1)$-th derivative of the first summand in \eqref{lt_H-(infty)-H-} is
given by formulae \eqref{Hxy0(infty)-Hxy0(lambda)_derivative} and \eqref{H0y0(infty)-H0y0(lambda)_derivative}.
Let us analyze the asymptotic properties of the $([d/2]-1)$-th derivative of the second summand in
\eqref{lt_H-(infty)-H-}. Integration by parts gives
\begin{eqnarray}\label{check_Hxy0_derivative}
\int\nolimits_{0}^{\infty}{e^{-\lambda t}t^{[d/2]-1}\,d H_{x,y,{\bf 0}}(t)}=&-&\lambda\int\nolimits_{0}^{\infty}
{e^{-\lambda t}t^{[d/2]-1}(H_{x,y,{\bf 0}}(\infty)-H_{x,y,{\bf 0}}(t))\,dt}\nonumber\\
&+&\left(\left[\frac{d}{2}\right]-1
\right)\int\nolimits_{0}^{\infty}{e^{-\lambda t}t^{[d/2]-2}(H_{x,y,{\bf 0}}(\infty)-H_{x,y,{\bf 0}}(t))\,dt}.\quad
\end{eqnarray}
Obviously, the first summand in \eqref{check_Hxy0_derivative} is $o\left(\int\nolimits_{0}^{\infty}
{e^{-\lambda t}t^{[d/2]-1}(H_{x,y,{\bf 0}}(\infty)-H_{x,y,{\bf 0}}(t))\,d t}\right)$, as ${\lambda\to0+}$.
By formulae \eqref{H_xyz(infty)-H_xyz(t)sim,d>=3}, \eqref{Hxy0(infty)-Hxy0(lambda)_derivative},
\eqref{H0y0(infty)-H0y0(lambda)_derivative} as well as Tauberian theorem
(Theorem 2 in \cite{Feller}, Ch.13, Sec.5) we infer that the second summand in \eqref{check_Hxy0_derivative}
is also $o\left(\int\nolimits_{0}^{\infty}{e^{-\lambda t}t^{[d/2]-1}(H_{x,y,{\bf 0}}(\infty)-
H_{x,y,{\bf 0}}(t))\,d t}\right)$,
as $\lambda\to0+$. Thus, the $([d/2]-1)$-th derivative of the second summand in \eqref{lt_H-(infty)-H-}
does not make a contribution to the asymptotic behavior of the $([d/2]-1)$-th derivative
of the left-hand side of \eqref{lt_H-(infty)-H-}. Therefore, $\int\nolimits_{0}^{\infty}
{e^{-\lambda t}t^{[d/2]-1}(H^{-}_{x,y,{\bf 0}}(\infty)-H^{-}_{x,y,{\bf 0}}(t))\,d t}\sim
\int\nolimits_{0}^{\infty}
{e^{-\lambda t}t^{[d/2]-1}(H_{x,y,{\bf 0}}(\infty)-H_{x,y,{\bf 0}}(t))\,d t}$ where the asymptotic relation
for the last expression is provided by formulae \eqref{Hxy0(infty)-Hxy0(lambda)_derivative} and
\eqref{H0y0(infty)-H0y0(lambda)_derivative}. Consequently, employing Corollary 43 in \cite{V_Lectures}
we complete the proof of Theorem \ref{T:H-} when $d\geq4$. Theorem \ref{T:H-} is proved.
\vskip0.2cm
The author is grateful to Associate Professor E.B. Yarovaya for permanent attention and to
Professor V.A. Vatutin for valuable remarks. Special thanks are to Professors I.Kourkova and
G.Pag\`es for invitation to LPMA UPMC.

\newpage
Ekaterina Vl. BULINSKAYA, \vskip0,2cm Faculty of Mathematics and
Mechanics,

Lomonosov Moscow State University,

Moscow 119991, Russia

\vskip0,5cm

{\it E-mail address}: bulinskaya@yandex.ru

\end{document}